\documentclass[11pt]{amsart}

\usepackage{pstricks}

\usepackage{amsfonts}
\usepackage{amsmath}
\usepackage{amsrefs}
\usepackage{amsthm}
\usepackage{array}
\usepackage{color}
\usepackage{graphics}
\usepackage{pstricks,pst-node}

\newcounter{movetype}

\newcounter{thm}

\newtheorem*{conjecture}{Conjecture}
\newtheorem{proposition}[thm]{Proposition}
\newtheorem{corollary}[thm]{Corollary}
\newtheorem{lemma}[thm]{Lemma}
\newtheorem{theorem}[thm]{Theorem}

\definecolor{grey}{gray}{.4}

\def\bbox#1{\framebox(#1,#1){$\bullet$}}

\def\ebox#1{\framebox(#1,#1){}}
\def\fbox#1{\framebox(#1,#1){\textcolor{grey}{$\bullet$}}}

\def\C{{\mathbb C}}                 
\def\Cc{{\mathcal C}}               
\def\dlow{\mathrm{dlow}}            
\def\dhigh{\mathrm{dhigh}}          
\def\F{{\mathbb F}}                 
\def\G{{\mathbb G}}                 
\def\l{{\ell}}                      
\def\lhat{{\hat{\ell}}}             
\def\L{{\mathcal L}}                
\def\Lhat{{\hat{\mathcal L}}}       
\def\O{{\mathcal O}}                
\def\P{{\mathbb P}}                 
\def\R{{\mathcal R}}                
\def\vhat{{\hat{v}}}                

\title{Linear systems in $\mathbb{P}^2$ with base points of bounded multiplicity}
\author{Stephanie Yang}

\begin{document}
\maketitle
\begin{abstract}
We present a proof of the Harbourne-Hirschowitz conjecture for linear
systems with multiple points of order 7 or less.  This uses a well-known
degeneration of the plane developed by Ciliberto and Miranda as well
as a combinatorial game that arises from specializing points onto
lines.
\end{abstract}



\section{Introduction}\label{sec:intro}

This paper discusses two techniques for determining if general
multiple points in $\P^2_\C$ impose independent linear conditions on
the space of plane curves of a given degree.  A well-known conjecture,
formulated independently by Harbourne, Gimigliano, and
Hirschowitz,~\cites{MR846019,MR91a:14007,MR993223} gives geometric
meaning to when this is the case.

Let $m_1,\ldots,m_r$ be a sequence of positive integers corresponding
to general points $p_1,\ldots,p_r\in\P^2$.  Denote by
$\L=\L_d(m_1,\ldots,m_k)$ the linear system of degree $d$ curves with
multiplicity $m_i$ at $p_i$.  Vanishing to order $m$ at a point $p$ is
equivalent to the vanishing of all derivatives at $p$ of order at most
$m-1$.  Thus an $m$-fold point imposes $\binom{m+1}{2}$ linear
conditions on plane curves and the expected dimension of $\L$ is given
by the equation:
\begin{align}
   e(\L)&=\max\left\{\binom{d+2}{2}-\sum_{i=1}^k\binom{m_i+1}{2}-1,-1\right\}.
\end{align}
This estimate is a sharp lower bound for the actual dimension of $\L$.
When equality holds, we say that $\L$ is {\em non-special}, and
otherwise, we say that $\L$ is {\em special}.

Let $\pi\colon V \to\P^2$ be the blow-up of the projective plane at the
points $p_1,\ldots,p_r$.  A curve $C\subseteq\P^2$ is called a {\em
$(-1)$-curve} if it is rational and its proper transform $\widetilde
C\subseteq V$ has self-intersection equal to $-1$.  With this in our
vocabulary, it is easy to state the Harbourne-Hirschowitz conjecture:

\begin{conjecture}[Harbourne-Hirschowitz]
$\L$ is special if and only if it contains a multiple (-1)-curve in
its base locus.
\end{conjecture}
\noindent While one direction (the ``if'' part) of this equivalence is
elementary, the other direction remains open except for special cases.

The Harbourne-Hirschowitz conjecture has a variety of
algebro-geometric consequences.  First, a proof of the conjecture
would settle the longstanding Nagata conjecture, posed in 1959 by
Nagata after he constructed a counterexample to Hilbert's 14th
problem~\cite{MR21:4151}.  (In short, Nagata conjectured if $n\geq
10$, then any degree $d$ curve with $n$ points of multiplicity $m$
must satisfy $d>m\sqrt{n}$).  The Harbourne-Hirschowitz conjecture
also implies that any integral curve of negative self-intersection in
the blow-up of $\P^2$ at (any number of) general points must have
self-intersection -1, thus giving a complete description of the Mori
cone of such surfaces.

One approach to this problem has a simple geometric description.
Suppose we are given a linear system $\L$ of plane curves with
multiple base points.  Choose a triangle of three lines in $\P^2$
meeting in three distinct points.  We specialize the base points by
moving them onto these points and sliding the multiple points along
the three lines to collide them.  Each collision creates a larger
singularity in the base locus of the limiting linear system, and the
class of singularities that arise can be completely described via a
combinatorial game involving checkers on a triangular board.

The second technique is a modification of a planar degeneration
exploited by Ciliberto and Miranda in~\cite{MR2000m:14005}
and~\cite{MR2000m:14006}.  Let $\Delta$ be a one-parameter family, and
denote by $X$ the blow-up of the three-fold $\P^2\times\Delta$ at a
point.  The fibers of $X$ over $\Delta$ can be viewed as a family of
projective planes which degenerate to a reducible surface comprised of
two rational components.  If we have a family of plane curves with
multiple points in $\P^2$, we can use this degeneration to `break' a
family of plane curves into two families defined on each of the two
rational components of the special fiber of $X$.  This gives a
recursive bound for the dimension of the original family.  A
consequence of this degeneration is the following statement, made
precise in Theorem~\ref{th:induction}.

\medskip
{\em
\centerline{
For any positive integer $M$, there exists a bound
$D=D(M)$ such that:
}
\smallskip
$\left.
\parbox{140pt}{The Harbourne-Hirschowitz conjecture is true for all
  linear systems $\L_d(m_1,\ldots,m_k)$ with $d<D(M)$ and $m_i\leq
  M$.}  \right\}
\Longrightarrow\left\{
\parbox{140pt}{The Harbourne-Hirschowitz conjecture is true for all
  linear systems $\L_d(m_1,\ldots,m_k)$ with $m_i\leq M$.}\right.$
}
\medskip

\noindent
The base points above are allowed to have mixed multiplicity.  (Most
recent results have applied only to collections of base points with
all or all but one points of equal multiplicity.)  Also note that the
list of linear systems for the left hand side is finite in length,
while those on the right are infinite.  The exact formula for $D(M)$
is given by Equation~\ref{eq:dm} in Section~\ref{sec:mixed}.

In particular, $D(7)=29$ and the number of possible linear systems
$28$ or less, with multiple points of order $7$ or less, is
approximately $10^8$.  One hundred million cases sounds daunting to
all but the computer-minded.  We wrote a program (in \texttt{C++}) to
enumerate this long list of cases and play the combinatorial game (of
checkers on a triangle) on each case.  Remarkably, the game worked to
prove the Harbourne-Hirschowitz conjecture in almost all of the cases,
cutting the number down to 42 (listed in Table~\ref{42}), which are
then handled with ad hoc methods in the last section of this paper, to
prove:

\begin{theorem}\label{th:seven}
The Harbourne-Hirschowitz conjecture is true for all linear systems of
plane curves with base points having multiplicity at most $7$.
\end{theorem}


\section{Checkers on a triangular board}\label{sec:game}

\subsection{Combinatorial description}

In this section, we describe the rules of a combinatorial game which
involves placing and moving up to $\sum \binom{m_i+1}{2}$ checkers on a
triangular checkerboard with side length $d+1$, containing a total of
$\binom{d+2}{2}$ squares.  The ultimate goal of the game is to place
as many checkers on the board as possible; this gives an upper bound
for the dimension of a linear system $\L_d(m_1,\ldots,m_r)$.

Given $\L_d(m_1,\ldots,m_k)$, form a $(d+1)\times (d+1)$ triangle of
boxes.  We may place checkers in these boxes using only two types of
moves:

\begin{list}{\textsc{Type \Alph{movetype}:}}{\usecounter{movetype}}

\item For any multiplicity $m_i$, we place $\binom{m_i+1}{2}$ checkers
in one of the three corners of the box, forming an $m_i\times m_i$
triangle.  If no corner of the box has enough empty squares available,
then our only options are to quit the game, or perform moves of type B
in order to create more empty squares in a corner.  Two examples of
valid moves are:

\bigskip
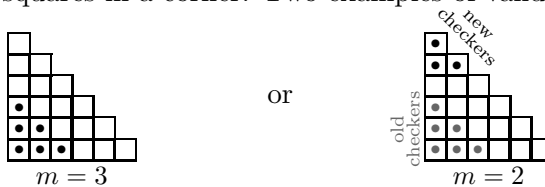
\hfill
\begin{picture}(48,48)\tiny
  \put(0,0){\bbox{8}}
  \put(0,8){\bbox{8}}
  \put(0,16){\bbox{8}}
  \put(0,24){\ebox{8}}
  \put(0,32){\ebox{8}}
  \put(0,40){\ebox{8}}
  \put(8,0){\bbox{8}}
  \put(8,8){\bbox{8}}
  \put(8,16){\ebox{8}}
  \put(8,24){\ebox{8}}
  \put(8,32){\ebox{8}}
  \put(16,0){\bbox{8}}
  \put(16,8){\ebox{8}}
  \put(16,16){\ebox{8}}
  \put(16,24){\ebox{8}}
  \put(24,0){\ebox{8}}
  \put(24,8){\ebox{8}}
  \put(24,16){\ebox{8}}
  \put(32,0){\ebox{8}}
  \put(32,8){\ebox{8}}
  \put(40,0){\ebox{8}}
  \put(24,-6){\makebox(0,0){\small $m=3$}}
\end{picture}
\hfill
\begin{picture}(0,48)
  \put(0,24){\makebox(0,0){or}}
\end{picture}
\hfill
\begin{picture}(48,48)\tiny
  \put(0,0){\fbox{8}}
  \put(0,8){\fbox{8}}
  \put(0,16){\fbox{8}}
  \put(0,24){\ebox{8}}
  \put(0,32){\bbox{8}}
  \put(0,40){\bbox{8}}
  \put(8,0){\fbox{8}}
  \put(8,8){\fbox{8}}
  \put(8,16){\ebox{8}}
  \put(8,24){\ebox{8}}
  \put(8,32){\bbox{8}}
  \put(16,0){\fbox{8}}
  \put(16,8){\ebox{8}}
  \put(16,16){\ebox{8}}
  \put(16,24){\ebox{8}}
  \put(24,0){\ebox{8}}
  \put(24,8){\ebox{8}}
  \put(24,16){\ebox{8}}
  \put(32,0){\ebox{8}}
  \put(32,8){\ebox{8}}
  \put(40,0){\ebox{8}}
  \put(24,46){\makebox(0,0){
      \rotatebox{-45}{\parbox[b]{29pt}{new}}}}
  \put(15,47){\makebox(0,0){
      \rotatebox{-45}{\parbox[b]{29pt}{checkers}}}}
  \put(-12,21){\makebox(0,0){
      \rotatebox{90}{\parbox[b]{29pt}{\textcolor{grey}{old}}}}}
  \put(-6,12){\makebox(0,0){
      \rotatebox{90}{\parbox[b]{29pt}{\textcolor{grey}{checkers}}}}}
 \put(24,-6){\makebox(0,0){\small $m=2$}}
\end{picture}
\hfill
\bigskip

\item We may perform one of six `slides' which move all of the
checkers as far as possible and in the same direction along rows,
columns, or diagonals.  The checkers may not overlap or share squares.
Two examples of valid moves of this type are:

\bigskip
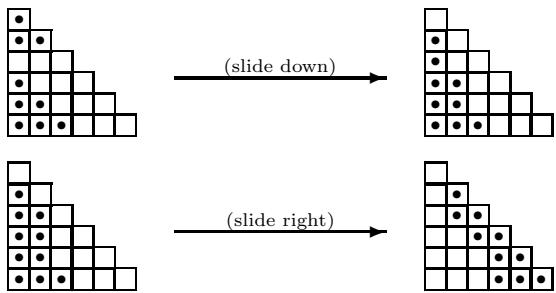
\hfill
\begin{picture}(48,48)\tiny
  \put(0,0){\bbox{8}}
  \put(0,8){\bbox{8}}
  \put(0,16){\bbox{8}}
  \put(0,24){\ebox{8}}
  \put(0,32){\bbox{8}}
  \put(0,40){\bbox{8}}
  \put(8,0){\bbox{8}}
  \put(8,8){\bbox{8}}
  \put(8,16){\ebox{8}}
  \put(8,24){\ebox{8}}
  \put(8,32){\bbox{8}}
  \put(16,0){\bbox{8}}
  \put(16,8){\ebox{8}}
  \put(16,16){\ebox{8}}
  \put(16,24){\ebox{8}}
  \put(24,0){\ebox{8}}
  \put(24,8){\ebox{8}}
  \put(24,16){\ebox{8}}
  \put(32,0){\ebox{8}}
  \put(32,8){\ebox{8}}
  \put(40,0){\ebox{8}}
\end{picture}
\hfill
\begin{picture}(0,48)\small\thicklines
  \put(0,26){\makebox(0,0){\tiny(slide down)}}
  \put(-40,22){\vector(1,0){80}}
\end{picture}
\hfill
\begin{picture}(48,48)\tiny
  \put(0,0){\bbox{8}}
  \put(0,8){\bbox{8}}
  \put(0,16){\bbox{8}}
  \put(0,24){\bbox{8}}
  \put(0,32){\bbox{8}}
  \put(0,40){\ebox{8}}
  \put(8,0){\bbox{8}}
  \put(8,8){\bbox{8}}
  \put(8,16){\bbox{8}}
  \put(8,24){\ebox{8}}
  \put(8,32){\ebox{8}}
  \put(16,0){\bbox{8}}
  \put(16,8){\ebox{8}}
  \put(16,16){\ebox{8}}
  \put(16,24){\ebox{8}}
  \put(24,0){\ebox{8}}
  \put(24,8){\ebox{8}}
  \put(24,16){\ebox{8}}
  \put(32,0){\ebox{8}}
  \put(32,8){\ebox{8}}
  \put(40,0){\ebox{8}}
\end{picture}
\hfill
\bigskip

\hfill
\begin{picture}(48,48)\tiny
  \put(0,0){\bbox{8}}
  \put(0,8){\bbox{8}}
  \put(0,16){\bbox{8}}
  \put(0,24){\bbox{8}}
  \put(0,32){\bbox{8}}
  \put(0,40){\ebox{8}}
  \put(8,0){\bbox{8}}
  \put(8,8){\bbox{8}}
  \put(8,16){\bbox{8}}
  \put(8,24){\bbox{8}}
  \put(8,32){\ebox{8}}
  \put(16,0){\bbox{8}}
  \put(16,8){\ebox{8}}
  \put(16,16){\ebox{8}}
  \put(16,24){\ebox{8}}
  \put(24,0){\ebox{8}}
  \put(24,8){\ebox{8}}
  \put(24,16){\ebox{8}}
  \put(32,0){\ebox{8}}
  \put(32,8){\ebox{8}}
  \put(40,0){\ebox{8}}
\end{picture}
\hfill
\begin{picture}(0,48)\small\thicklines
  \put(0,26){\makebox(0,0){\tiny(slide right)}}
  \put(-40,22){\vector(1,0){80}}
\end{picture}
\hfill
\begin{picture}(48,48)\tiny
  \put(0,0){\ebox{8}}
  \put(0,8){\ebox{8}}
  \put(0,16){\ebox{8}}
  \put(0,24){\ebox{8}}
  \put(0,32){\ebox{8}}
  \put(0,40){\ebox{8}}
  \put(8,0){\ebox{8}}
  \put(8,8){\ebox{8}}
  \put(8,16){\ebox{8}}
  \put(8,24){\bbox{8}}
  \put(8,32){\bbox{8}}
  \put(16,0){\ebox{8}}
  \put(16,8){\ebox{8}}
  \put(16,16){\bbox{8}}
  \put(16,24){\bbox{8}}
  \put(24,0){\bbox{8}}
  \put(24,8){\bbox{8}}
  \put(24,16){\bbox{8}}
  \put(32,0){\bbox{8}}
  \put(32,8){\bbox{8}}
  \put(40,0){\bbox{8}}
\end{picture}
\hfill
\bigskip
\end{list}
The dimension of $\L_d(m_1,\ldots,m_r)$ is bounded above by one less
than the number of uncheckered boxes after any sequence of moves, so
long as we use each multiplicity $m_i$ for moves of type A at most
once.  In other words, if all of the $\sum\binom{m_i+1}{2}$ checkers
can be fit into the triangle using only the two moves described above,
then $\L_d(m_1,\ldots,m_r)$ is non-special.

\medskip
\noindent
{\bf Examples.}
As a first example, consider linear system $\L_5(3,2,2,2,2,2)$ of quintics
with one triple point and five double points.  When we perform the
triangle algorithm as follows:

\begin{list}{}{}
\item[\textsc{Step 1:}] Place six checkers (for the triple point) onto the lower right
hand corner of boxes.

\item[\textsc{Steps 2--3:}] Place three checkers for a double
point onto the lower left hand corner of boxes, and slide all the
checkers to the right.

\item[\textsc{Steps 4--5:}] Place three checkers for a double
point onto the upper corner of boxes, and slide all the checkers down.

\item[\textsc{Steps 6--8:}] Place three checkers for a double point on
the top corner, slide the checkers down, and then slide them to the right.

\item[\textsc{Steps 9--10}] Repeat steps 4 and 5.

\item[\textsc{Step 11:}] Place three more checkers, for the last
double point, into the remaining empty boxes in the upper corner.
\end{list}

\noindent On the triangular checkerboard, the steps look like this:

\bigskip
\noindent
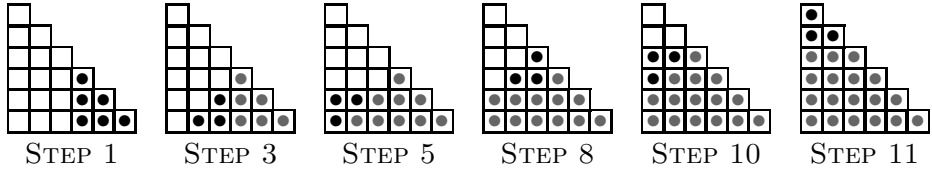
\hfill
\begin{picture}(48,56)(0,-8)
   \put(0,0){\ebox{8}}
   \put(8,0){\ebox{8}}
   \put(16,0){\ebox{8}}
   \put(24,0){\bbox{8}}
   \put(32,0){\bbox{8}}
   \put(40,0){\bbox{8}}
   \put(0,8){\ebox{8}}
   \put(8,8){\ebox{8}}
   \put(16,8){\ebox{8}}
   \put(24,8){\bbox{8}}
   \put(32,8){\bbox{8}}
   \put(0,16){\ebox{8}}
   \put(8,16){\ebox{8}}
   \put(16,16){\ebox{8}}
   \put(24,16){\bbox{8}}
   \put(0,24){\ebox{8}}
   \put(8,24){\ebox{8}}
   \put(16,24){\ebox{8}}
   \put(0,32){\ebox{8}}
   \put(8,32){\ebox{8}}
   \put(0,40){\ebox{8}}
   \put(24,-8){\makebox(0,0){\textsc{Step 1}}}
\end{picture}
\hfill
\begin{picture}(48,56)(0,-8)
   \put(0,0){\ebox{8}}
   \put(8,0){\bbox{8}}
   \put(16,0){\bbox{8}}
   \put(24,0){\fbox{8}}
   \put(32,0){\fbox{8}}
   \put(40,0){\fbox{8}}
   \put(0,8){\ebox{8}}
   \put(8,8){\ebox{8}}
   \put(16,8){\bbox{8}}
   \put(24,8){\fbox{8}}
   \put(32,8){\fbox{8}}
   \put(0,16){\ebox{8}}
   \put(8,16){\ebox{8}}
   \put(16,16){\ebox{8}}
   \put(24,16){\fbox{8}}
   \put(0,24){\ebox{8}}
   \put(8,24){\ebox{8}}
   \put(16,24){\ebox{8}}
   \put(0,32){\ebox{8}}
   \put(8,32){\ebox{8}}
   \put(0,40){\ebox{8}}
   \put(24,-8){\makebox(0,0){\textsc{Step 3}}}
\end{picture}
\hfill
\begin{picture}(48,56)(0,-8)
   \put(0,0){\bbox{8}}
   \put(8,0){\fbox{8}}
   \put(16,0){\fbox{8}}
   \put(24,0){\fbox{8}}
   \put(32,0){\fbox{8}}
   \put(40,0){\fbox{8}}
   \put(0,8){\bbox{8}}
   \put(8,8){\bbox{8}}
   \put(16,8){\fbox{8}}
   \put(24,8){\fbox{8}}
   \put(32,8){\fbox{8}}
   \put(0,16){\ebox{8}}
   \put(8,16){\ebox{8}}
   \put(16,16){\ebox{8}}
   \put(24,16){\fbox{8}}
   \put(0,24){\ebox{8}}
   \put(8,24){\ebox{8}}
   \put(16,24){\ebox{8}}
   \put(0,32){\ebox{8}}
   \put(8,32){\ebox{8}}
   \put(0,40){\ebox{8}}
   \put(24,-8){\makebox(0,0){\textsc{Step 5}}}
\end{picture}
\hfill
\begin{picture}(48,56)(0,-8)
   \put(0,0){\fbox{8}}
   \put(8,0){\fbox{8}}
   \put(16,0){\fbox{8}}
   \put(24,0){\fbox{8}}
   \put(32,0){\fbox{8}}
   \put(40,0){\fbox{8}}
   \put(0,8){\fbox{8}}
   \put(8,8){\fbox{8}}
   \put(16,8){\fbox{8}}
   \put(24,8){\fbox{8}}
   \put(32,8){\fbox{8}}
   \put(0,16){\ebox{8}}
   \put(8,16){\bbox{8}}
   \put(16,16){\bbox{8}}
   \put(24,16){\fbox{8}}
   \put(0,24){\ebox{8}}
   \put(8,24){\ebox{8}}
   \put(16,24){\bbox{8}}
   \put(0,32){\ebox{8}}
   \put(8,32){\ebox{8}}
   \put(0,40){\ebox{8}}
   \put(24,-8){\makebox(0,0){\textsc{Step 8}}}
\end{picture}
\hfill
\begin{picture}(48,56)(0,-8)
   \put(0,0){\fbox{8}}
   \put(8,0){\fbox{8}}
   \put(16,0){\fbox{8}}
   \put(24,0){\fbox{8}}
   \put(32,0){\fbox{8}}
   \put(40,0){\fbox{8}}
   \put(0,8){\fbox{8}}
   \put(8,8){\fbox{8}}
   \put(16,8){\fbox{8}}
   \put(24,8){\fbox{8}}
   \put(32,8){\fbox{8}}
   \put(0,16){\bbox{8}}
   \put(8,16){\fbox{8}}
   \put(16,16){\fbox{8}}
   \put(24,16){\fbox{8}}
   \put(0,24){\bbox{8}}
   \put(8,24){\bbox{8}}
   \put(16,24){\fbox{8}}
   \put(0,32){\ebox{8}}
   \put(8,32){\ebox{8}}
   \put(0,40){\ebox{8}}
   \put(24,-8){\makebox(0,0){\textsc{Step 10}}}
\end{picture}
\hfill
\begin{picture}(48,56)(0,-8)
   \put(0,0){\fbox{8}}
   \put(8,0){\fbox{8}}
   \put(16,0){\fbox{8}}
   \put(24,0){\fbox{8}}
   \put(32,0){\fbox{8}}
   \put(40,0){\fbox{8}}
   \put(0,8){\fbox{8}}
   \put(8,8){\fbox{8}}
   \put(16,8){\fbox{8}}
   \put(24,8){\fbox{8}}
   \put(32,8){\fbox{8}}
   \put(0,16){\fbox{8}}
   \put(8,16){\fbox{8}}
   \put(16,16){\fbox{8}}
   \put(24,16){\fbox{8}}
   \put(0,24){\fbox{8}}
   \put(8,24){\fbox{8}}
   \put(16,24){\fbox{8}}
   \put(0,32){\bbox{8}}
   \put(8,32){\bbox{8}}
   \put(0,40){\bbox{8}}
   \put(24,-8){\makebox(0,0){\textsc{Step 11}}}
\end{picture}
\hfill

\bigskip
\noindent The darker dots $\bullet$ represent the newly placed
checkers, while the lighter dots $\textcolor{grey}\bullet$ represent
checkers from previous moves.  All 21 checkers fit into the board, and
thus $\L_5(3,2,2,2,2,2)$ is empty and non-special.

Now consider the special linear system $\L=\L_2(2,2)$ of conics
through two double points.  After placing first three checkers in any
corner of the triangle, we cannot fit another triangle of three
checkers onto the board, even after any sequence of slides.  This of
course is due to the fact that $\L_2(2,2)$ is special.

\medskip
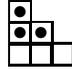
\hfill
\begin{picture}(24,24)
   \put(0,0){\ebox{8}}
   \put(0,8){\bbox{8}}
   \put(0,16){\bbox{8}}
   \put(8,0){\ebox{8}}
   \put(8,8){\bbox{8}}
   \put(16,0){\ebox{8}}
\end{picture}
\hfill
\medskip

\subsection{Algebraic description of the game}

We translate these moves into the language of algebra using pictures
to guide us with bookkeeping.  Choose homogeneous coordinates for
$\P^2$ so that the vertices of the triangle are $[1\colon 0\colon 0]$,
$[0\colon 1\colon 0]$, $[0\colon 0\colon 1]$.  Using these
coordinates, we can create a $(d+1)\times(d+1)$ triangle of boxes
which represent the monomial basis for degree $d$ curves, with $X^d$,
$Y^d$, $Z^d$ represented by the three corner boxes, and the other
monomials interpolated between these in the usual way.

\medskip
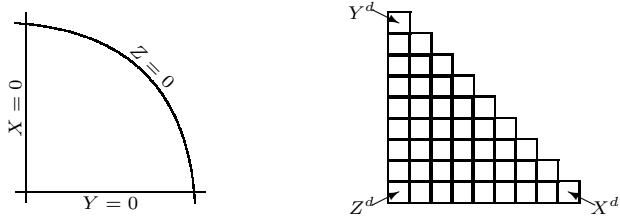
\hfill
\begin{picture}(72,72)\tiny
  \put(4,0){\line(0,1){72}}
  \put(0,4){\line(1,0){72}}
  \qbezier(0,68)(64,64)(68,0)
  \put(36,0){\makebox(0,0){$Y=0$}}
  \put(0,36){\rotatebox{90}{\makebox(0,0){{$X=0$}}}}
  \put(52,52){\rotatebox{-45}{\makebox(0,0){{$Z=0$}}}}
\end{picture}
\hfill
\begin{picture}(72,72)\tiny
  \put(0,0){\ebox{8}}
  \put(0,8){\ebox{8}}
  \put(0,16){\ebox{8}}
  \put(0,24){\ebox{8}}
  \put(0,32){\ebox{8}}
  \put(0,40){\ebox{8}}
  \put(0,48){\ebox{8}}
  \put(0,56){\ebox{8}}
  \put(0,64){\ebox{8}}
  \put(8,0){\ebox{8}}
  \put(8,8){\ebox{8}}
  \put(8,16){\ebox{8}}
  \put(8,24){\ebox{8}}
  \put(8,32){\ebox{8}}
  \put(8,40){\ebox{8}}
  \put(8,48){\ebox{8}}
  \put(8,56){\ebox{8}}
  \put(16,0){\ebox{8}}
  \put(16,8){\ebox{8}}
  \put(16,16){\ebox{8}}
  \put(16,24){\ebox{8}}
  \put(16,32){\ebox{8}}
  \put(16,40){\ebox{8}}
  \put(16,48){\ebox{8}}
  \put(24,0){\ebox{8}}
  \put(24,8){\ebox{8}}
  \put(24,16){\ebox{8}}
  \put(24,24){\ebox{8}}
  \put(24,32){\ebox{8}}
  \put(24,40){\ebox{8}}
  \put(32,0){\ebox{8}}
  \put(32,8){\ebox{8}}
  \put(32,16){\ebox{8}}
  \put(32,24){\ebox{8}}
  \put(32,32){\ebox{8}}
  \put(40,0){\ebox{8}}
  \put(40,8){\ebox{8}}
  \put(40,16){\ebox{8}}
  \put(40,24){\ebox{8}}
  \put(48,0){\ebox{8}}
  \put(48,8){\ebox{8}}
  \put(48,16){\ebox{8}}
  \put(56,0){\ebox{8}}
  \put(56,8){\ebox{8}}
  \put(64,0){\ebox{8}}

  \put(-10,-1){\makebox(0,0){$Z^d$}}
  \put(-6,-1){\vector(2,1){10}}
  \put(-10,73){\makebox(0,0){$Y^d$}}
  \put(-6,73){\vector(2,-1){10}}
  \put(82,-1){\makebox(0,0){$X^d$}}
  \put(78,-1){\vector(-2,1){10}}
\end{picture}
\hfill
\medskip

\noindent
Notice that imposing an $m$-tuple base point onto $[1\colon 0\colon
0]$, $[0\colon 1\colon 0]$, or $[0\colon 0\colon 1]$ corresponds
exactly to the vanishing of the monomials in the an $m\times m$ corner
of the triangle.  We denote this by placing a checker in the boxes
that represent the vanishing monomials:

\medskip
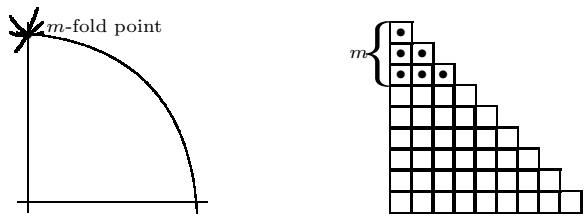
\hfill
\begin{picture}(72,72)\tiny
  \put(4,0){\line(0,1){72}}
  \put(0,4){\line(1,0){72}}
  \qbezier(0,68)(64,64)(68,0)
  \thicklines
  \qbezier(-2.5,68)(3.5,65)(9.5,72)
  \qbezier(-.5,61)(3.5,65)(7.5,77)
  \qbezier(-2.5,74)(3.5,65)(9.5,68)
  \put(4,67.3){\circle*{3}}
  \put(33,70){\makebox(0,0){$m$-fold point}}
\end{picture}
\hfill
\begin{picture}(72,72)\tiny
  \put(0,0){\ebox{8}}
  \put(0,8){\ebox{8}}
  \put(0,16){\ebox{8}}
  \put(0,24){\ebox{8}}
  \put(0,32){\ebox{8}}
  \put(0,40){\ebox{8}}
  \put(0,48){\bbox{8}}
  \put(0,56){\bbox{8}}
  \put(0,64){\bbox{8}}
  \put(8,0){\ebox{8}}
  \put(8,8){\ebox{8}}
  \put(8,16){\ebox{8}}
  \put(8,24){\ebox{8}}
  \put(8,32){\ebox{8}}
  \put(8,40){\ebox{8}}
  \put(8,48){\bbox{8}}
  \put(8,56){\bbox{8}}
  \put(16,0){\ebox{8}}
  \put(16,8){\ebox{8}}
  \put(16,16){\ebox{8}}
  \put(16,24){\ebox{8}}
  \put(16,32){\ebox{8}}
  \put(16,40){\ebox{8}}
  \put(16,48){\bbox{8}}
  \put(24,0){\ebox{8}}
  \put(24,8){\ebox{8}}
  \put(24,16){\ebox{8}}
  \put(24,24){\ebox{8}}
  \put(24,32){\ebox{8}}
  \put(24,40){\ebox{8}}
  \put(32,0){\ebox{8}}
  \put(32,8){\ebox{8}}
  \put(32,16){\ebox{8}}
  \put(32,24){\ebox{8}}
  \put(32,32){\ebox{8}}
  \put(40,0){\ebox{8}}
  \put(40,8){\ebox{8}}
  \put(40,16){\ebox{8}}
  \put(40,24){\ebox{8}}
  \put(48,0){\ebox{8}}
  \put(48,8){\ebox{8}}
  \put(48,16){\ebox{8}}
  \put(56,0){\ebox{8}}
  \put(56,8){\ebox{8}}
  \put(64,0){\ebox{8}}
  \put(-5,60){\makebox(0,0){\Huge\{}}
  \put(-12,60){\makebox(0,0){$m$}}
\end{picture}
\hfill
\medskip

Consider the six maps below; from now on we will refer to them as {\em
slide transformations}:
\begin{align*}
(X,Y,Z)&\mapsto (X+t^{-1}Y,Y,Z)
&(X,Y,Z)&\mapsto (X+t^{-1}Z,Y,Z)\\
(X,Y,Z)&\mapsto (X,Y+t^{-1}X,Z)
&(X,Y,Z)&\mapsto (X,Y+t^{-1}Z,Z)\\
(X,Y,Z)&\mapsto (X,Y,Z+t^{-1}X)
&(X,Y,Z)&\mapsto (X,Y,Z+t^{-1}Y)
\end{align*}
Let $\alpha$ denote a subset of $\{(i,j)\colon i+j\leq d\}$ and let
$I_\alpha$ denote the ideal generated by the set
$\{X^iY^jZ^{d-i-j}\colon (i,j)\in\alpha\}$.  Pictorially, we represent
$\alpha$ on our triangle by adding checkers to the monomials boxes
{\em not} in $I_\alpha$.  By lemma~\ref{le:slide}, the flat limit of
$I_\alpha$ under a slide transformation as $t$ vanishes is another
monomial ideal; we wish to calculate this flat limit and record how
this moves the checkers on our board.

Consider the slide transformation $(X,Y,Z)\mapsto (X+t^{-1}Y,Y,Z)$;
since this preserves the $Z$-degree of any monomial, we can simply
calculate how the slide transformation acts on the $Z$-homogeneous
parts of $I_\alpha$.  For this, assume that the set $\{(i,j)\colon
i+j\leq d\}-\alpha$ consists of pairs $(i,j)$ for which $i+j=N$ for
some fixed integer $N\leq d$.  In other words, we require all of the
checkered boxes on the triangle associated with $\alpha$ lie on a
diagonal, for example:

\medskip
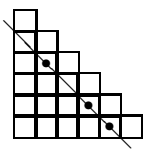
\hfill
\begin{picture}(48,48)\tiny
  \put(0,0){\ebox{8}}
  \put(0,8){\ebox{8}}
  \put(0,16){\ebox{8}}
  \put(0,24){\ebox{8}}
  \put(0,32){\ebox{8}}
  \put(0,40){\ebox{8}}
  \put(8,0){\ebox{8}}
  \put(8,8){\ebox{8}}
  \put(8,16){\ebox{8}}
  \put(8,24){\bbox{8}}
  \put(8,32){\ebox{8}}
  \put(16,0){\ebox{8}}
  \put(16,8){\ebox{8}}
  \put(16,16){\ebox{8}}
  \put(16,24){\ebox{8}}
  \put(24,0){\ebox{8}}
  \put(24,8){\bbox{8}}
  \put(24,16){\ebox{8}}
  \put(32,0){\bbox{8}}
  \put(32,8){\ebox{8}}
  \put(40,0){\ebox{8}}
  \put(44,-4){\line(-1,1){48}}
\end{picture}
\hfill
\medskip

\noindent
Lemma~\ref{le:slide} states that the flat limit of $I_\alpha$ under
the slide transformation $(X,Y,Z)\mapsto(X+t^{-1}Y,Y,Z)$ as $t$
vanishes is the ideal with all of the checkers boxes shifted along the
$Z$-homogeneous lines towards the highest power of $Y^d$.

\medskip
\hfill
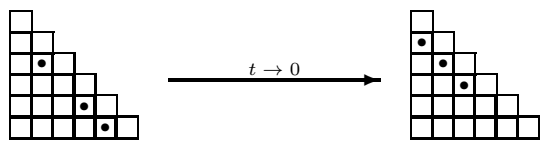
\begin{picture}(48,48)\tiny
  \put(0,0){\ebox{8}}
  \put(0,8){\ebox{8}}
  \put(0,16){\ebox{8}}
  \put(0,24){\ebox{8}}
  \put(0,32){\ebox{8}}
  \put(0,40){\ebox{8}}
  \put(8,0){\ebox{8}}
  \put(8,8){\ebox{8}}
  \put(8,16){\ebox{8}}
  \put(8,24){\bbox{8}}
  \put(8,32){\ebox{8}}
  \put(16,0){\ebox{8}}
  \put(16,8){\ebox{8}}
  \put(16,16){\ebox{8}}
  \put(16,24){\ebox{8}}
  \put(24,0){\ebox{8}}
  \put(24,8){\bbox{8}}
  \put(24,16){\ebox{8}}
  \put(32,0){\bbox{8}}
  \put(32,8){\ebox{8}}
  \put(40,0){\ebox{8}}
\end{picture}
\hfill
\begin{picture}(0,48)\small\thicklines
  \put(0,26){\makebox(0,0){\tiny$t\to 0$}}
  \put(-40,22){\vector(1,0){80}}
\end{picture}
\hfill
\begin{picture}(48,48)\tiny
  \put(0,0){\ebox{8}}
  \put(0,8){\ebox{8}}
  \put(0,16){\ebox{8}}
  \put(0,24){\ebox{8}}
  \put(0,32){\bbox{8}}
  \put(0,40){\ebox{8}}
  \put(8,0){\ebox{8}}
  \put(8,8){\ebox{8}}
  \put(8,16){\ebox{8}}
  \put(8,24){\bbox{8}}
  \put(8,32){\ebox{8}}
  \put(16,0){\ebox{8}}
  \put(16,8){\ebox{8}}
  \put(16,16){\bbox{8}}
  \put(16,24){\ebox{8}}
  \put(24,0){\ebox{8}}
  \put(24,8){\ebox{8}}
  \put(24,16){\ebox{8}}
  \put(32,0){\ebox{8}}
  \put(32,8){\ebox{8}}
  \put(40,0){\ebox{8}}
\end{picture}
\hfill\hfill
\medskip

This analysis is similar for all six slide transformations with the
appropriate permutation of $X$, $Y$, $Z$.  As a corollary, we have
that the flat limit if any ideal $I_\alpha$ under any slide
transformation results in sliding the checkers along lines parallel to
the sides of the triangle.

\begin{lemma}\label{le:slide}
Suppose the ideal $I_\alpha$ is generated by $l$ homogeneous monomials
of degree $k$ in their $Z$-degree, for some fixed $k$.  Then the flat
limit of the ideal $I_\alpha$ under the slide transformation
$T_t:(X,Y)\mapsto(X+t^{-1}Y,Y,Z)$ as $t$ vanishes is the monomial
ideal $I_\beta$ where $\beta=\{(i,d-k-i,k):i=0,\ldots,l-1\}$
\end{lemma}
\begin{proof}
Assume that $\alpha$ consists of pairs $(a_i,d-k-a_i)$ for
$i=1,\ldots,l$, with $a_1<a_2<\ldots,a_l$.  The $T_tI_\alpha$ is
generated by the $l$ elements:
\begin{equation}
   (X+t^{-1}Y)^{a_i}Y^{d-k-a_i}Z^k=\sum_{j=0}^{a_i}
   \binom{a_i}{j} t^{j-a_i} X^j Y^{d-k-j}Z^k
\end{equation}
Thinking of the set $\left\{X^iY^{d-k-i}Z^k\right\}_{i=0}^d$ as a
basis for the vector space of degree $d$ monomials in $X,Y,Z$ with
$Z$-degree equal to $k$, we can form $l\times d$ matrix representing
the equations above:
\begin{equation}
M=\begin{bmatrix}
   \binom{a_1}{0}t^{-a_1} & \binom{a_1}{1}t^{1-a_1}
   & \cdots & \binom{a_0}{d}t^{d-a_0} &\\
   \binom{a_2}{0}t^{-a_2} & \binom{a_2}{1}t^{1-a_2}
   & \cdots & \binom{a_1}{d}t^{d-a_1} &\\
   \vdots & \vdots & \ddots & \vdots \\
   \binom{a_l}{0}t^{-a_l} & \binom{a_l}{1}t^{1-a_l} & \ldots &
   \binom{a_l}{d}t^{d-a_l} &
\end{bmatrix}
\end{equation}
The first $k$ columns of this matrix are linearly independent.  To see
this, we compute the determinant
\begin{align}
\begin{vmatrix}
   \binom{a_0}{0}t^{-a_0} & \binom{a_0}{1}t^{1-a_0}
   & \cdots & \binom{a_0}{l}t^{l-a_0} &\\
   \binom{a_1}{0}t^{-a_1} & \binom{a_1}{1}t^{1-a_1}
   & \cdots & \binom{a_1}{l}t^{l-a_1} &\\ \vdots & \vdots &
   \ddots & \vdots \\ \binom{a_l}{0}t^{a_l} &
   \binom{-a_l}{1}t^{1-a_l} & \ldots &
   \binom{-a_l}{l}t^{l-a_l} &
\end{vmatrix} &= 
t^{-N}
\begin{vmatrix}
  \binom{a_0}{0} & \binom{a_0}{1} & \cdots &
  \binom{a_l}{l} &\\
  \binom{a_1}{0} & \binom{a_1}{1} & \cdots & \binom{a_1}{l} &\\ \vdots
  & \vdots & \ddots & \vdots \\
  \binom{a_l}{0} & \binom{a_l}{1} & \ldots &
  \binom{a_l}{l} &
\end{vmatrix} \\
&= \frac{1}{1!\cdots l!}\prod_{i<j}(a_i-a_j)
\end{align}
where $N=\sum{a_i}-\frac{k(k+1)}{2}\geq 0$.  The last equality comes
from the fact that we can factor the second matrix above as an
lower-diagonal matrix times a Van der Monde matrix:
\begin{equation}
\begin{bmatrix}
   \binom{a_1}{0} &  \binom{a_2}{0} & \cdots & \binom{a_l}{0} &\\
   \binom{a_1}{1} &  \binom{a_2}{1} & \cdots & \binom{a_l}{1} &\\
   \vdots & \vdots & \ddots & \vdots \\
   \binom{a_1}{l} &  \binom{a_2}{l} & \ldots & \binom{a_l}{l} &
\end{bmatrix}
=
\begin{bmatrix}
1 & 0 & 0 & 0 & \cdots \\
0 & 1 & 0 & 0 & \cdots \\
0 & -\frac{1}{2} & \frac{1}{2} &0 & \cdots \\
0 & \frac{1}{3}  &-\frac{1}{2} &\frac{1}{6} &\cdots \\
\vdots &\vdots &\vdots &\vdots &\ddots
\end{bmatrix}
\begin{bmatrix}
   1 & 1 & 1 & \cdots & 1 \\
   a_1 & a_2 & a_3 &\cdots &a_l \\
   a_1^2 & a_2^2 &a_3^2 &\cdots &a_l^2 \\
   a_1^3 & a_2^3 &a_3^3 &\cdots &a_l^3 \\
   \vdots & \vdots &\ddots &\vdots \\
   a_1^l & a_2^l &a_3^l &\cdots &a_l^l
\end{bmatrix}
\end{equation}
Thus we are able to perform row operations on the matrix $M$ to obtain
the identity matrix in the first $l$ rows.  If this is done, all of
the elements of $M$ outside the first $k$ columns will contain a power
of $t$.  In other words, we can re-express the basis of
$T_tI_\alpha$:
\begin{align}
 T_t I_\alpha &= \langle
 X+t^{-1}Y)^{a_i}Y^{d-k-a_i}Z^k\rangle_{i=1}^l \\ &= \langle
 X^iY^{d-k-i}Z^k+\text{\small higher order terms in $t$}\rangle_{i=0}^{l-1}
\end{align}
and the flat limit of this ideal as $t$ vanishes is simply
$X^iY^{d-k-i}Z^k$ as $i=0,\ldots,l-1$.
\end{proof}

By upper-semicontinuity, the dimension of $\L_d(m_1,\ldots,m_k)$ is
bounded above by one less than the dimension of the monomial ideal
which arises from playing the triangular checkers game with all of
them multiplicities $m_i$; this is exactly the number of white boxes
left at the end of the game.  Thus, if there is some order of moves
which fits all of the checkers onto the triangular board, then
$\L_d(m_1,\ldots,m_k)$ is non-special.


\section{Degenerating $\P^2$}\label{sec:cm}

\subsection{Degenerating the plane}\label{sec:blowup}

In this section, we describe a degeneration of $\P^2$ used by
Ciliberto and Miranda in~\cites{MR2000m:14005,MR2000m:14006} and then
by them with Cioffi and Orecchia in~\cite{CCMO} to prove the
Harbourne-Hirschowitz conjecture for homogeneous linear systems with
$m\leq 20$.  The degeneration described here is only slightly modified
from their degeneration for the purpose of including linear systems
with base points of mixed multiplicity.

Let $\Delta$ be a disc around the origin, and let $\pi:
X\to\Delta\times\P^2$ be the three-fold obtained by blowing-up the
product $\Delta\times\P^2$ at a point $p$ in the plane
$\{0\}\times\P^2$.

Denote by $X_t=\pi^{-1}(t)$ the fiber of $X$ over $t\in\Delta$.  If
$t\neq 0$, then $X_t \cong \P^2$.  The special fiber $X_0$ is the
union of the exceptional divisor $\P$, which is a copy of the
projective plane, with the Hirzebruch space $\F$, isomorphic to $\P^2$
blown-up at a point.  Via this isomorphism it is easy to see that the
Picard group of $\F$ is freely generated by two divisors $H$ and $E$,
where $H$ is the pullback of a class of a general line in $\P^2$ and
$E$ is the class of the exceptional divisor of the blow-up.  Let
$R\subseteq X_0$ denote the divisor $\P\cap\F$.

\begin{equation}
\begin{psmatrix}[colsep=.15cm,rowsep=.8cm]
\P\cup_R\F=&X_0 &&&X                      \\
               &&&&\Delta\times\P^2       \\
&0              &&&\Delta           &&&\P^2
\psset{arrows=->,nodesep=3pt}
\ncline{1,2}{1,5}
\ncline{1,2}{3,2}
\ncline{1,5}{2,5}
\ncarc[arcangle=20,labelsep=0]{1,5}{3,8}>{\phi}
\ncline{2,5}{3,8}
\ncline{2,5}{3,5}
\ncline{3,2}{3,5}
\end{psmatrix}
\end{equation}

Let $\phi:X \to \P^2$ denote the composition of $\pi$ with the
standard projection from the second factor of $\Delta\times\P^2$.
Denote by $\O(d,a)$ the line bundle
\begin{align}
   \O(d,a)=\phi^*\O_{\P^2}(d)\otimes\O_X(-a\P).
\end{align}
For any $0\leq a \leq d$, the line bundle $\O(d,a)$ is a flat family
of line bundles over $\Delta$.  If $t\neq 0$, then this line bundle
restricts to $X_t$ as $\O_{\P^2}(d)$.  If $t=0$, then this line bundle
restricts to the components of the special fiber $X_0$ as
\begin{align}
   \O(d,a)|_{\P}  &= \O_{\P}(a) \\
   \O(d,a)|_{\F}  &= \O_{\F}(dH-aE),
\end{align}
This follows from the fact that $\O_X(\P)$ restricts to $\P$ as
$\O_{\P^2}(-1)$ and to $\F$ as $\O_\F(E)$.  These two line bundles
above agree when restricted to $R$.

We now modify our notation to help us index collections of multiple
points in $\P^2$.  Let $\L_d(m_1^{k_1},\ldots,m_s^{k_s})$ denote the
linear system of degree $d$ curves through $\sum k_i$ general points,
of which $k_i$ have multiplicity $m_i$.  (For example,
$\L_7(3,3,3,3,3,3)$ can be abbreviated to $\L_7(3^6)$, and
$\L_5(2,2,2,3)$ can be written as $\L_5(2^3,3)$.)  Let
$l_1,\ldots,l_s$ be another sequence of positive integers such that
$l_i \leq k_i$ for $i=1,\ldots,s$.

Consider $\sum k_i$ general points in the reducible fiber $X_0$ with
$l_i$ of the $m_i$-fold points in $\F$, for $i=1,\ldots,s$, and the
rest of the points in $\P$, all in general position.  These points can
be considered as limits of a family of multiple points in general
position in the nearby fibers of $\pi$.  Denote by $\L_0$ the linear
system of divisors in $|\O(d,a)|$ which vanish at these multiple
points in $X_0$.  By semicontinuity we have
\begin{align}
   \dim\L_0\geq \dim \L_d(m_1^{k_1},\ldots,m_s^{k_s}).
\end{align}
Our goal is to find parameters $a$ and $l_i$ that make $\dim \L_0$ is
as small as possible.

We take advantage of the fact that $X_0$ is a reducible surface whose
components are rational.  Specifically, let $\L_{\P}$ and $\L_{\F}$
denote the restrictions of $\L_0$ to $\P$ and $\F$.  Then,
\begin{align}
   \L_{\P}&\cong\L_a(m_1^{k_1-l_1},\ldots,m_s^{k_s-l_s}) \\
   \L_{\F}&\cong\L_d(m_1^{l_1},\ldots,m_s^{l_s},a).
\end{align}
The second equation comes from blowing down the $(-1)$-curve $E$ in
$\F$.  Of course, $a$ can be equal to any of the multiplicities $m_i$;
for convenience of notation it is easier to leave it as is.  Let
$\R_\P$ and $\R_\F$ denote the restrictions of $\L_\P$ and $\L_\F$ to
$R$, and denote by $\Lhat_\P$ and $\Lhat_\F$ denote the kernels of
these restrictions.  Then
\begin{align}
   \Lhat_{\P}&\cong\L_{a-1}(m_1^{k_1-l_1},\ldots,m_s^{k_s-l_s}) \\
   \Lhat_{\F}&\cong\L_d(m_1^{l_1},\ldots,m_s^{l_s},a+1).
\end{align}
Denote by $\l_\P$, $\l_\F$, $\lhat_\P$, $\lhat_\F$, $r_\P$, and $r_\F$
the dimension of the respective linear systems $\L_\P$, $\L_\F$,
$\Lhat_\P$, $\Lhat_\F$, $\R_\P$, and $\R_\F$.  A study of these linear
systems in~\cite{MR2000m:14006} offers the following equations:
\begin{align}
   r_\P &= \l_\P - \lhat_\P - 1 \\
   r_\F &= \l_\F - \lhat_\F - 1 \\
   \l_0 &= \dim(\R_\P\cap\R_\F)+\lhat_\P+\lhat_\F+1. \label{eq:diml0}
\end{align}
In the last equation, $\dim(\R_\P\cap\R_\F)$ refers to the vector
space dimension of the intersection of two linear systems inside the
$(a+1)$-dimensional $H^0\O_R(a)$.  The first two
equations are immediate from definitions.  The idea behind
equation~\ref{eq:diml0} is that elements of $\L_0$ come from pairs of
elements in $\L_\P$ and $\L_\F$ which agree on $R$.  If $\L_\P = \P
W_\P$ and $\L_\F=\P W_\F$, then
\begin{align}
\L_0 = \P \left(W_\P \times_{H^0 \O_R(a)} W_\F\right),
\end{align}
and the result follows from a simple dimension count.

If the systems $\R_\P$ and $\R_\F$ intersect transversely, then
$\dim(\R_\P\cap\R_\F)$ is immediate, and we would have a recursion for
$\l_0$.  In fact, transverse intersection is always the case, thanks
to a proof by Hirschowitz in~\cite{MR932136} and again by Ciliberto
and Miranda in~\cite{MR2000m:14005}.  This gives us the following
proposition and corollary:
\begin{proposition}\label{pro:l0}
With the notation as above,
\begin{enumerate}
\item If $r_\P+r_\F\leq a+1$, then $\l_0=\lhat_\P+\lhat_\F+1$.
\item If $r_\P+r_\F> a+1$, then $\l_0=\l_\P+\l_\F-a$.
\end{enumerate}
\end{proposition}

\begin{corollary}\label{co:l0}
If there exists $a$ and $l_i$ such that, with the notation above,
$\L_\P$ and $\L_\F$ are non-special, and $\Lhat_\P$ and $\Lhat_\F$ are
empty, then $\L_d(m_1^{k_1},\ldots,m_s^{k_s})$ is non-special.
\end{corollary}

\begin{proof}
We can assume that $e(\L_d(m_1^{k_1},\ldots,m_s^{k_s})) \leq 0$,
otherwise, we can impose additional simple base points until this is
the case.  If $\lhat_\P=\lhat_\F=-1$, then that $r_\P=\l_\P$ and
$r_\F=\l_\F$.  Thus
\begin{align}
\begin{split}
   r_\P+r_\F=&
      \binom{a+2}{2}-\sum k_i\binom{m_i+1}{2}+\sum l_i\binom{m_i+1}{2} \\
      &+\binom{d+2}{2}-\sum l_i\binom{m_i+1}{2}-\binom{a+1}{2}
\end{split}\\
      \leq&\ a+1 + e\left(\L_d(m_1^{k_1},\ldots,m_s^{k_s})\right).
\end{align}
The condition for first part of Proposition~\ref{pro:l0} is satisfied,
and consequently $\l_0=\lhat_\P+\lhat_\F+1=-1$ as expected.
\end{proof}

We end this section with a final useful observation, which will be
used several times in section~\ref{sec:main}.  Let $\Delta^* =
\Delta-\{0\}$, and suppose $\L_d(m_1^{k_1},\ldots,m_s^{k_s})$ contains
a curve which lies in some fiber $X_t\cong \P^2$ for some
$t\in\Delta^*$.  As we vary the base points of this curve, it sweeps
out a surface in $\Cc^*\subseteq\P^2\times\Delta^*$.  The closure of
$\pi^{-1}C^*\subseteq X$ we denote by $\Cc$.

\begin{lemma}\label{le:match}
Let $B$ be a (-1)-curve in $X_0$ such that $B\cdot R = 1$.  If
$B\cdot\Cc = -\sigma < 0$, then $\Cc$ contains $B$ to order at least
$\sigma$.
\end{lemma}

The details of this proof, as well as generalizations of the statement
above, are discussed in~\cite{MR2129797} and~\cite{matching}.  The
idea behind the proof is to blow up $X$ along the $r$ sections above
$\Delta$ which are the general points in each $\P^2$, and then to blow
up again along the proper transform of $B$.  The exceptional divisor
$\G$ of the last blow-up is a quadric surface; a quick calculation
expressing the restriction to $\G$ of the proper transform of $\Cc$
yields the result.

These ``matching conditions'' are the basis of a useful technique that
relates the non-speciality of a linear system
$\L_d(m_1^{k_1},\ldots,m_s^{k_s})$ to the non-speciality of a linear
system of lower degree with points in special position.  An example of
this technique is carried out in Lemma~\ref{le:match} of
Section~\ref{sec:main}.

\subsection{The induction argument for large degree}\label{sec:mixed}

Define the following rather unsightly functions,
\begin{align}
   \dlow(\gamma,h,m)  &= \left\lceil
                            \frac{\binom{m}{2}+\binom{\gamma+1}{2}
                            +(2h+1)\binom{m+1}{2}-m\gamma-1}
                            {m+1-\gamma}
                         \right\rceil \\
   \dhigh(\gamma,h,m) &= m+h+mh+\gamma h-1\\
   S(M)               &= \left\lceil-\frac{3}{2} +
                            \left(\sum_{m=1}^M \left(2\left\lceil
                               \frac{m^2-1}{3m+4}
                            \right\rceil+1\right)m(m+1)\right)^{\frac{1}{2}}
                         \right\rceil.
\end{align}
In this section we prove the following theorem.

\begin{theorem}\label{th:induction}
Let $M = \max\{m_i\}$, and let
\begin{align}\label{eq:dm}
  D = D(M) = \max\left\{4M+1,
         \dlow\left(-1,\left\lceil\frac{M^2-1}{3M+4}\right\rceil,M\right),
              S(M) \right\}.
\end{align}
Suppose the Harbourne-Hirschowitz conjecture holds for all linear systems
$\L_d(m_1^{k_1},\ldots,m_s^{k_s})$ with $m_i\leq M$ for $i=1,\ldots,s$
and $d<D(M)$.  Then it is true for all
$\L_d(m_1^{k_1},\ldots,m_s^{k_s})$ with $m_i\leq M$ for $i=1,\ldots,s$
and all values of $d$.
\end{theorem}
This is based on an induction on both the degree $d$ and the
multiplicities $m_i$ of the points, starting with the fact that if and
$d\geq 3m_i$ for all $i$, then $\L$ contains no (-1)-curves.  Together
with the Harbourne-Hirschowitz conjecture, this fact implies that any
special linear system must have at least one base point of
multiplicity $d/3$ or greater.

\begin{lemma}\label{le:3m}
If $\L_d(m_1,\ldots,m_r)$ contains a $(-1)$-curve in its base locus,
then $d<3\max\{m_i\}$.
\end{lemma}

\begin{proof}
Let $M=\max\{m_i\}$, and let $C$ denote a $(-1)$-curve in the base
locus of $\L$.  As before, let $V$ be the blow-up of $\P^2$ at the
base points $p_i$ of $\L$, and suppose that $\tilde C$, the proper
transform of $C$, represents the class $eE_0-\sum n_iE_i$, where $E_0$
denotes the class of the pullback of a general line and $E_i$ denotes
the class of the exceptional divisor over $p_i$.  The adjunction
formula gives us the equation
\begin{align}
   3e -\sum n_i &= 1.
\end{align}
If $|D|$ contains $\tilde C$, then $\tilde C \cdot D < 0$, and thus
\begin{align}
   de &< \sum n_i m_i \leq M \sum n_i = M (3e-1), \\
   d  &< 3M.
\end{align}
\end{proof}

A linear system is called {\em quasi-homogeneous} if all of its base
points except one are of the same multiplicity; i.e, if
it is of the form $\L_d(m^k,a)$.  A study of quasi-homogeneous systems
in~\cite{MR2000m:14005} yields the following result.

\begin{lemma}\label{le:quasi}
$\L_d(m^b,d-m+\gamma)$ is non-special if, $2 \leq m \leq d$, $b$ is
odd, and $-1\leq\gamma\leq 1$.
\end{lemma}

The proof of this lemma follows from analysis of the equations that
arise from the definition of quasi-homogeneous (-1)-curves and those
that arise from permuting the multiple points of the same order within
in a quasi-homogeneous systems.  The lemma for $\gamma=0$ and
$\gamma=1$ can also be seen easily using the methods discussed in
section~\ref{sec:game}; the details are left to the reader.

\medskip

\noindent
{\em Proof of Theorem~\ref{th:induction}.}  As in the proof of
Corollary~\ref{co:l0}, we may assume that
$e\left(\L_d(m_1^{k_1},\ldots,m_s^{k_s})\right)\leq 0$; otherwise, we add
general simple points until this is the case.  Assume the
Harbourne-Hirschowitz conjecture is true for all
$\L_d(m_1^{k_1},\ldots,m_s^{k_s})$ with $m_i< M$ and $d<D$.  We will
find parameters $a$ and $l_i$ so that the systems $\L_\P$, $\L_\F$,
$\Lhat_\P$ and $\Lhat_\F$ satisfy the conditions of
Corollary~\ref{co:l0}.

Let $l_i=2h+1$ be odd and $l_j=0$ for $j\neq i$, and set
$a=d-m_i-\gamma$, for $\gamma \in \{-1,0,1\}$.  By Lemma~\ref{le:quasi},
the linear system $\L_\F=\L_d(m_i^{l_i},a)$ is non-special, and the
linear system $\Lhat_\F=\L_d(m_i^{l_i},a+1)$ contains $\gamma+1$ lines
through the $(a+1)$-fold point and the $m_i$-fold points.  The residual
system to these lines is
\begin{align}
   \L_{d-l_i(\gamma+1)}\left((m_i-\gamma-1)^{l_i},a+1-l_i(\gamma+1)\right),
\end{align}
which is also non-special by Lemma~\ref{le:quasi}.

If $d\geq D$, then $a\geq 3M$, and thus $\L_\P$ and $\Lhat_\P$ are
also non-special by the induction hypothesis and Lemma~\ref{le:3m}.

Thus all four linear systems $\L_\P$, $\L_\F$, $\Lhat_\P$, and
$\Lhat_\F$ are non-special, and their virtual dimensions are:
\begin{align}
   v_\P &=\binom{a+2}{2}-\sum_{j=1}^s l_j\binom{m_j+1}{2}
         +l_i\binom{m_i+1}{2}-1 \\
   v_\F &=\binom{d+2}{2}-\binom{a+1}{2}-l_i\binom{m_i+1}{2}-1 \\
   \vhat_\P &=\binom{a+1}{2}-\sum_{j=1}^s l_j\binom{m_j+1}{2}
         +l_i\binom{m_i+1}{2}-1 \\
   \vhat_\F &=\binom{d+2}{2}-\binom{a+2}{2}-l_i\binom{m_i+1}{2}-1.
\end{align}
For $\Lhat_\P$ and $\Lhat_\F$ to be empty, we need that $d$ (and
consequently $a$) be small enough that $\vhat_\P\leq 0$, yet large
enough that $\vhat_\F\leq 0$.  Imposing $\vhat_\P\leq 0$ and
$\vhat_\F\leq 0$ yield, respectively, an upper and lower bounds for
$d$:
\begin{align}
   d&\leq\dlow(\gamma,h,m_i) \\
   d&\geq \dhigh(\gamma,h,m_i).
\end{align}
Thus for fixed $a=d-m_i+\gamma$ and $l_i=2h+1$, we have an interval
\begin{align}
   \dlow(\gamma,h,m_i) \leq d \leq \dhigh(\gamma,h,m_i)
\end{align}
on which we can prove that $\L_d(m_1^{k_1},\ldots,m_s^{k_s})$ is
empty. We force these intervals to overlap by varying $\gamma$ and
$h$.

It is a ``remarkable fact'' \cite{MR2000m:14005} that
$\dhigh(-1,h,m_i)=\dlow(0,h,m_i)$ and 
$\dhigh(0,h,m_i) =\dlow(1,h,m_i)$.
This extends our induction interval to
\begin{align}
   \dlow(m_i,-1,h) \leq d \leq \dhigh(m_i,1,h).
\end{align}
We now vary $h$ so that these intervals overlap, which gives us the
equations
\begin{align}
   \dhigh(m_i,1,h) + 1 &\geq \dlow(m_i,-1,h) \\
   h &\geq \left\lceil \frac{m_i^2-1}{3m_i+4} \right\rceil.
\end{align}
For this to work, we must have the existence of some $i$ for which
$k_i \geq 2h+1 = 2 \left\lceil\frac{m_i^2-1}{3m_i+4}\right\rceil+1$.  Since
$\chi\left(\L_d(m_1^{k_1},\ldots,m_s^{k_s})\right)\leq 0$, this last requirement
is fulfilled for large $d$.
\begin{align}
   \binom{d+2}{2}
      &\geq \sum_{i=1}^M\left(2\left\lceil
                           \frac{i^2+i}{3i+4}
                        \right\rceil+1\right)
                                \binom{i+1}{2} \\
    d &\geq S(M).
\end{align}
\hfill {$\Box$}


\section{Proof of Theorem~\ref{th:seven}}\label{sec:main}

Below, Table~\ref{dm} shows the first few values of $D(M)$.  For all
but the first few values of $M$, the function $D(M)$ is determined by
$S(M)$.
\begin{table}[ht]
\hfill
\begin{tabular}{c||c|c|c||c}
$M$ & $4M+1$
   &$\dlow\left(-1,\left\lceil\frac{M^2-1}{3M+4}\right\rceil,M\right)$
   &$S(M)$
   &$D(M)$ \\ \hline
 2 & 9 & 3 & 3 &  9 \\
 3 &13 & 5 & 6 & 13 \\
 4 &17 & 7 &10 & 17 \\
 5 &21 &13 &15 & 21 \\
 6 &25 &16 &21 & 25 \\
 7 &29 &19 &26 & 29 \\
 8 &33 &29 &34 & 34 \\
 9 &37 &33 &42 & 42 \\
10 &41 &37 &51 & 51 \\
11 &45 &51 &61 & 61 \\
12 &49 &56 &71 & 71
\end{tabular}
\hfill\hfill
\caption{Values of $D(M)$ for $M=2,\ldots,12$}\label{dm}
\end{table}
To prove Theorem~\ref{th:seven}, we programmed a computer to enumerate
all linear systems $\L_d(m_1,\ldots,m_k)$ of degree 29 or less, with
points of multiplicity $7$ or less.  There $125220076$ of these,
almost all of which were shown to satisfy the Harbourne-Hirschowitz
conjecture via the game discussed in Section~\ref{sec:game}:
\begin{center}
\begin{tabular}{cr>{\small}l}
  & 125220076 &  total linear systems\\
--& 124850912 & are empty via the combinatorial method\\
--& 366691 & are empty because otherwise there would \\
  &        & \ \ appear ``too many curves'' in the base locus\\
--& 2013 & are special because of multiple (-1)-curves in the base locus\\
--& 418 & are empty using Proposition~\ref{pro:l0}\\\hline
 & 42 & systems remain and are listed on Table~\ref{42}
\end{tabular}
\end{center}

\begin{table}
\begin{center}
\begin{tabular}{>{\small}l|>{\small}l|>{\small}l}
 $\L_d(m_1^{k_1},\ldots,m_s^{k_s})$ & Reason for being empty & Order
 of specialization \\ \hline
 $\L_{15}(3, 4, 5^{8})$ & Triangular checker game & $5^5, 4, 3$ \\ 
 $\L_{16}(2, 5^{10})$ & Triangular checker game & $5^5, 2$, then $4^2$\\ 
 $\L_{17}(1^{2}, 5^{8}, 6, 7)$ & Cremona transformation \\ 
 $\L_{17}(2, 5^{7}, 6^{3})$ & Cremona transformation \\ 
 $\L_{17}(4, 5^{7}, 7^{2})$ & Cremona transformation \\ 
 $\L_{18}(1, 3, 5, 6^{8})$ & Triangular checker game &$6^5, 5, 3, 1$ \\ 
 $\L_{18}(1, 5^{7}, 7^{3})$ & Cremona transformation \\ 
 $\L_{18}(1^{2}, 5^{6}, 6^{2}, 7^{2})$ & Cremona transformation \\ 
 $\L_{18}(2, 6^{9})$ & Cubics in the base locus \\ 
 $\L_{18}(4^{5}, 7^{5})$ & Cremona transformation \\ 
 $\L_{19}(1, 3, 6^{7}, 7^{2})$ & Cremona transformation \\ 
 $\L_{19}(1, 6^{10})$ & Homogeneous \\ 
 $\L_{19}(2^{2}, 5, 6^{9})$ & Triangular checker game & $5,6^6,2^2$ \\ 
 $\L_{19}(3, 5, 6, 7^{6})$ & Cremona transformation \\ 
 $\L_{19}(3, 5, 6^{5}, 7^{3})$ & Cremona transformation \\ 
 $\L_{19}(3, 5, 6^{9})$ & Triangular checker game &$6^6, 5, 3$ \\ 
 $\L_{19}(4, 5^{4}, 7^{5})$ & Cremona transformation \\ 
 $\L_{19}(5, 6^{8}, 7)$ & Triangular checker game &$7, 6^4, 5, 6$ \\ 
 $\L_{19}(6^{10})$ & Homogeneous \\ 
 $\L_{20}(1, 3, 6^{4}, 7^{5})$ & Cremona transformation \\ 
 $\L_{20}(1, 3, 6^{8}, 7^{2})$ & Triangular checker game &$7^2, 6^5, 1, 3$ \\ 
 $\L_{20}(1, 6^{11})$ & Homogeneous \\ 
 $\L_{20}(1, 6^{7}, 7^{3})$ & Cremona transformation \\ 
 $\L_{20}(3, 5, 6^{2}, 7^{6})$ & Cremona transformation \\ 
 $\L_{20}(3, 5, 6^{6}, 7^{3})$ & Cremona transformation \\ 
 $\L_{20}(5, 6^{5}, 7^{4})$ & Cremona transformation \\ 
 $\L_{20}(6^{11})$ & Homogeneous \\ 
 $\L_{20}(6^{7}, 7^{3})$ & Cremona transformation \\ 
 $\L_{21}(1, 2, 4, 5, 7^{8})$ & Triangular checker game &$7^4, 1, 7, 5, 2$ \\ 
 $\L_{21}(1, 6^{4}, 7^{6})$ & Lemma~\ref{le:two} \\ 
 $\L_{21}(1^{2}, 3, 6, 7^{8})$ & Lemma~\ref{le:two} \\ 
 $\L_{21}(2, 7^{9})$ & Cubics in base locus \\ 
 $\L_{21}(5, 6^{2}, 7^{7})$ & Lemma~\ref{le:two} \\ 
 $\L_{22}(1, 2, 6, 7^{9})$ & Implied by $\L_{22}(2, 6, 7^{9})$ \\ 
 $\L_{22}(1^{3}, 6, 7^{9})$ & Implied by $\L_{22}(2, 6, 7^{9})$  \\ 
 $\L_{22}(2, 6, 7^{9})$ & Lemma~\ref{le:two} \\ 
 $\L_{22}(2, 6^{13})$ & Lemma~\ref{le:save} \\ 
 $\L_{22}(4, 5, 7^{9})$ & Triangular checker game & $5,7^5,4$ \\ 
 $\L_{22}(4, 6^{2}, 7^{8})$ & Triangular checker game & $7^3,6,7^2,6,4$ \\ 
 $\L_{23}(6, 7^{10})$ & Triangular checker game & $7^7, 6$ \\ 
 $\L_{26}(6^{18})$ & Homogeneous \\ 
 $\L_{27}(5, 7^{14})$ & Triangular checker game & $7^{11}, 5$ \\ 
\end{tabular}
\end{center}
\caption{The forty two linear systems leftover from the computer program}\label{42}
\end{table}

Twelve of the linear systems in the Table~\ref{42} can be handled via
a combination of the triangular checker game and an analysis of which
curves are forced to appear in the base locus.  For example, to show that $\L_{15}(3,4,5^8)$ contains
no curves, we play the triangle game in the following manner:
\begin{list}{\textsc{}}{}
  \item[\textsc{Steps 1--15:}] Place fifteen checkers (for a 5-tuple point)
  on the top of the triangle, and slide them down, and then slide them
  again to the right.  Repeat this four more times.

  \item[\textsc{Steps 16--18:}] Place ten checkers (for a quadruple
  point) on top of the triangle, slide them down, and then slide them
  to the right.

  \item[\textsc{Steps 19--21:}] Place six checkers (for a triple point)
  on top of the triangle, and slide them down, and then slide them to
  the right.
\end{list}
After these seven steps, a line $L$ splits off six times (since we
have six full rows of checkers on the bottom of the triangle).
Residual to the line, we are left with the monomial ideal represented
by the triangle below, and three unspecialized 5-tuple points:

\medskip
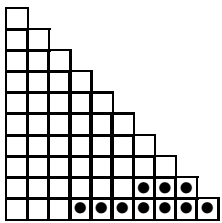
\hfill
\begin{picture}(80,80)
   \put(0,0){\ebox{8}}
   \put(8,0){\ebox{8}}
   \put(16,0){\ebox{8}}
   \put(24,0){\bbox{8}}
   \put(32,0){\bbox{8}}
   \put(40,0){\bbox{8}}
   \put(48,0){\bbox{8}}
   \put(56,0){\bbox{8}}
   \put(64,0){\bbox{8}}
   \put(72,0){\bbox{8}}
   \put(0,8){\ebox{8}}
   \put(8,8){\ebox{8}}
   \put(16,8){\ebox{8}}
   \put(24,8){\ebox{8}}
   \put(32,8){\ebox{8}}
   \put(40,8){\ebox{8}}
   \put(48,8){\bbox{8}}
   \put(56,8){\bbox{8}}
   \put(64,8){\bbox{8}}
   \put(0,16){\ebox{8}}
   \put(8,16){\ebox{8}}
   \put(16,16){\ebox{8}}
   \put(24,16){\ebox{8}}
   \put(32,16){\ebox{8}}
   \put(40,16){\ebox{8}}
   \put(48,16){\ebox{8}}
   \put(56,16){\ebox{8}}
   \put(0,24){\ebox{8}}
   \put(8,24){\ebox{8}}
   \put(16,24){\ebox{8}}
   \put(24,24){\ebox{8}}
   \put(32,24){\ebox{8}}
   \put(40,24){\ebox{8}}
   \put(48,24){\ebox{8}}
   \put(00,32){\ebox{8}}
   \put(8,32){\ebox{8}}
   \put(16,32){\ebox{8}}
   \put(24,32){\ebox{8}}
   \put(32,32){\ebox{8}}
   \put(40,32){\ebox{8}}
   \put(0,40){\ebox{8}}
   \put(8,40){\ebox{8}}
   \put(16,40){\ebox{8}}
   \put(24,40){\ebox{8}}
   \put(32,40){\ebox{8}}
   \put(0,48){\ebox{8}}
   \put(8,48){\ebox{8}}
   \put(16,48){\ebox{8}}
   \put(24,48){\ebox{8}}
   \put(0,56){\ebox{8}}
   \put(8,56){\ebox{8}}
   \put(16,56){\ebox{8}}
   \put(0,64){\ebox{8}}
   \put(8,64){\ebox{8}}
   \put(0,72){\ebox{8}}
\end{picture}
\hfill
\medskip

\noindent
Any degree $9$ curve containing the three unspecialized $5$-tuple
points must contain the lines through any two of the points; thus the
base locus of the new linear system contains three lines, as well as
the line $L$ once again.  Residual to these four lines, we have the
linear system of degree $5$ curves through three general triple
points, with a tangency to $L$ of order $3$.

\medskip
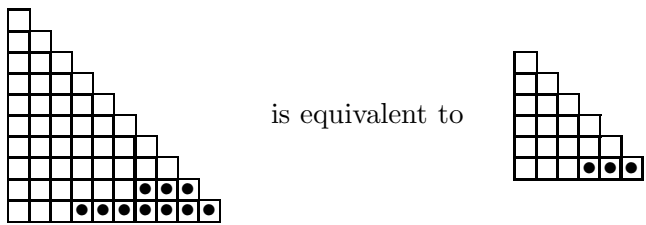
\hfill
\begin{picture}(80,80)
   \put(0,0){\ebox{8}}
   \put(8,0){\ebox{8}}
   \put(16,0){\ebox{8}}
   \put(24,0){\bbox{8}}
   \put(32,0){\bbox{8}}
   \put(40,0){\bbox{8}}
   \put(48,0){\bbox{8}}
   \put(56,0){\bbox{8}}
   \put(64,0){\bbox{8}}
   \put(72,0){\bbox{8}}
   \put(0,8){\ebox{8}}
   \put(8,8){\ebox{8}}
   \put(16,8){\ebox{8}}
   \put(24,8){\ebox{8}}
   \put(32,8){\ebox{8}}
   \put(40,8){\ebox{8}}
   \put(48,8){\bbox{8}}
   \put(56,8){\bbox{8}}
   \put(64,8){\bbox{8}}
   \put(0,16){\ebox{8}}
   \put(8,16){\ebox{8}}
   \put(16,16){\ebox{8}}
   \put(24,16){\ebox{8}}
   \put(32,16){\ebox{8}}
   \put(40,16){\ebox{8}}
   \put(48,16){\ebox{8}}
   \put(56,16){\ebox{8}}
   \put(0,24){\ebox{8}}
   \put(8,24){\ebox{8}}
   \put(16,24){\ebox{8}}
   \put(24,24){\ebox{8}}
   \put(32,24){\ebox{8}}
   \put(40,24){\ebox{8}}
   \put(48,24){\ebox{8}}
   \put(00,32){\ebox{8}}
   \put(8,32){\ebox{8}}
   \put(16,32){\ebox{8}}
   \put(24,32){\ebox{8}}
   \put(32,32){\ebox{8}}
   \put(40,32){\ebox{8}}
   \put(0,40){\ebox{8}}
   \put(8,40){\ebox{8}}
   \put(16,40){\ebox{8}}
   \put(24,40){\ebox{8}}
   \put(32,40){\ebox{8}}
   \put(0,48){\ebox{8}}
   \put(8,48){\ebox{8}}
   \put(16,48){\ebox{8}}
   \put(24,48){\ebox{8}}
   \put(0,56){\ebox{8}}
   \put(8,56){\ebox{8}}
   \put(16,56){\ebox{8}}
   \put(0,64){\ebox{8}}
   \put(8,64){\ebox{8}}
   \put(0,72){\ebox{8}}
\end{picture}
\hfill
\begin{picture}(0,80)
\put(0,40){\makebox(0,0){is equivalent to}}
\end{picture}
\hfill
\begin{picture}(48,80)(0,-16)
   \put(0,0){\ebox{8}}
   \put(8,0){\ebox{8}}
   \put(16,0){\ebox{8}}
   \put(24,0){\bbox{8}}
   \put(32,0){\bbox{8}}
   \put(40,0){\bbox{8}}
   \put(0,8){\ebox{8}}
   \put(8,8){\ebox{8}}
   \put(16,8){\ebox{8}}
   \put(24,8){\ebox{8}}
   \put(32,8){\ebox{8}}
   \put(0,16){\ebox{8}}
   \put(8,16){\ebox{8}}
   \put(16,16){\ebox{8}}
   \put(24,16){\ebox{8}}
   \put(0,24){\ebox{8}}
   \put(8,24){\ebox{8}}
   \put(16,24){\ebox{8}}
   \put(00,32){\ebox{8}}
   \put(8,32){\ebox{8}}
   \put(0,40){\ebox{8}}
\end{picture}
\hfill
\medskip

\noindent Once again, the triangle of lines through the three general
triple points and the line $L$ split off; the residual linear system
is simply $\L_1(3)$ which is clearly empty.

\medskip
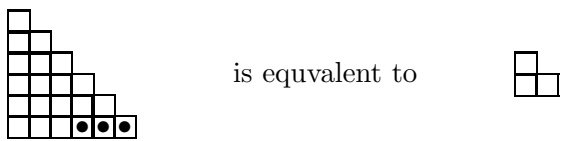
\hfill
\begin{picture}(48,48)
   \put(0,0){\ebox{8}}
   \put(8,0){\ebox{8}}
   \put(16,0){\ebox{8}}
   \put(24,0){\bbox{8}}
   \put(32,0){\bbox{8}}
   \put(40,0){\bbox{8}}
   \put(0,8){\ebox{8}}
   \put(8,8){\ebox{8}}
   \put(16,8){\ebox{8}}
   \put(24,8){\ebox{8}}
   \put(32,8){\ebox{8}}
   \put(0,16){\ebox{8}}
   \put(8,16){\ebox{8}}
   \put(16,16){\ebox{8}}
   \put(24,16){\ebox{8}}
   \put(0,24){\ebox{8}}
   \put(8,24){\ebox{8}}
   \put(16,24){\ebox{8}}
   \put(00,32){\ebox{8}}
   \put(8,32){\ebox{8}}
   \put(0,40){\ebox{8}}
\end{picture}
\hfill
\begin{picture}(0,48)
\put(0,24){\makebox(0,0){is equvalent to}}
\end{picture}
\hfill
\begin{picture}(16,48)(0,-16)
   \put(0,0){\ebox{8}}
   \put(8,0){\ebox{8}}
   \put(0,8){\ebox{8}}
\end{picture}
\hfill
\medskip

A similar argument works for eleven other linear systems in the table
which are empty via the triangular checker game; after specializing
all but three of the multiple points in the order prescribed by the
third column of Table~\ref{42}, triples of lines begin to split off of
the base locus.  Only $\L_{16}(2,5^{10})$ is slightly different.  In
this case, after specializing five quintuple points and a double
point, a conic appears through the five remaining quintuple points.
Residual to that, we have the ideal

\medskip
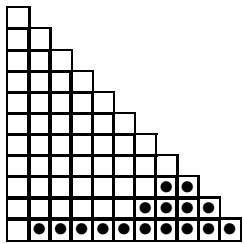
\hfill
\begin{picture}(88,88)
   \put(0,0){\ebox{8}}
   \put(8,0){\bbox{8}}
   \put(16,0){\bbox{8}}
   \put(24,0){\bbox{8}}
   \put(32,0){\bbox{8}}
   \put(40,0){\bbox{8}}
   \put(48,0){\bbox{8}}
   \put(56,0){\bbox{8}}
   \put(64,0){\bbox{8}}
   \put(72,0){\bbox{8}}
   \put(80,0){\bbox{8}}
   \put(0,8){\ebox{8}}
   \put(8,8){\ebox{8}}
   \put(16,8){\ebox{8}}
   \put(24,8){\ebox{8}}
   \put(32,8){\ebox{8}}
   \put(40,8){\ebox{8}}
   \put(48,8){\bbox{8}}
   \put(56,8){\bbox{8}}
   \put(64,8){\bbox{8}}
   \put(72,8){\bbox{8}}
   \put(0,16){\ebox{8}}
   \put(8,16){\ebox{8}}
   \put(16,16){\ebox{8}}
   \put(24,16){\ebox{8}}
   \put(32,16){\ebox{8}}
   \put(40,16){\ebox{8}}
   \put(48,16){\ebox{8}}
   \put(56,16){\bbox{8}}
   \put(64,16){\bbox{8}}
   \put(0,24){\ebox{8}}
   \put(8,24){\ebox{8}}
   \put(16,24){\ebox{8}}
   \put(24,24){\ebox{8}}
   \put(32,24){\ebox{8}}
   \put(40,24){\ebox{8}}
   \put(48,24){\ebox{8}}
   \put(56,24){\ebox{8}}
   \put(0,32){\ebox{8}}
   \put(8,32){\ebox{8}}
   \put(16,32){\ebox{8}}
   \put(24,32){\ebox{8}}
   \put(32,32){\ebox{8}}
   \put(40,32){\ebox{8}}
   \put(48,32){\ebox{8}}
   \put(0,40){\ebox{8}}
   \put(8,40){\ebox{8}}
   \put(16,40){\ebox{8}}
   \put(24,40){\ebox{8}}
   \put(32,40){\ebox{8}}
   \put(40,40){\ebox{8}}
   \put(0,48){\ebox{8}}
   \put(8,48){\ebox{8}}
   \put(16,48){\ebox{8}}
   \put(24,48){\ebox{8}}
   \put(32,48){\ebox{8}}
   \put(0,56){\ebox{8}}
   \put(8,56){\ebox{8}}
   \put(16,56){\ebox{8}}
   \put(24,56){\ebox{8}}
   \put(0,64){\ebox{8}}
   \put(8,64){\ebox{8}}
   \put(16,64){\ebox{8}}
   \put(0,72){\ebox{8}}
   \put(8,72){\ebox{8}}
   \put(0,80){\ebox{8}}
\end{picture}
\hfill
\medskip

\noindent
with five unspecialized quadruple points.  We specialize two of these
quadruple points, after triangles start to appear in the base locus as
in the previous argument.

Sixteen of the linear systems in the table are subject to quadratic
Cremona transformations; that is, suppose $m_1p_1$, $m_2p_2$, and
$m_3p_3$ are the three points of highest multiplicity.  We blow up
$\P^2$ at the three points $p_1$, $p_2$, and $p_3$ and then blow down
the resulting surface along the proper transforms of the lines
$\overline{p_1p_2}$, $\overline{p_1p_3}$, and $\overline{p_2p_3}$.  If
$s=m_1+m_2+m_3-d$, then the linear system $\L_d(m_1,\ldots,m_k)$ to
transformed to the linear system
$\L_{d-s}(m_1-s,m_2-s,m_3-s,m_4,\ldots,m_k)$.  In particular, if
$s>0$, then the Harbourne-Hirschowitz conjecture for
$\L_d(m_1,\ldots,m_k)$ is equivalent to the conjecture for a linear
system of lower degree (see~\cite{MR91a:14007}).

The last two lemmas of this section will introduce two more techniques
to prove the non-speciality of linear systems.  This will handle the
rest of the non-homogeneous cases in the table.

\begin{lemma}\label{le:two}
The linear systems $\L_{21}(1, 6^{4}, 7^{6})$, $\L_{21}(1^{2}, 3, 6, 7^{8})$, $\L_{21}(5,
 6^{2}, 7^{7})$, and $\L_{22}(2, 6, 7^{9})$ are empty.
\end{lemma}

\begin{proof}
The first three cases are so similar that we will only prove that one
of the linear systems is non-special; the rest follow by an almost
identical argument.  To show that $\L_{21}(1,6^4,7^6)$ is empty, we
degenerate $\P^2$ into a reducible surface $\P \cup \F$ as described
in Section~\ref{sec:cm}, placing four of the $7$-tuple points in $\F$
and the rest of the points in $\P$ and setting $a=17$.  Five curves
appear with multiplicity three in the base locus of $\L_\F$ (if
$\pi\colon\F\to\P^2$ is the blow-down of $\F$ along its $(-1)$-curve
$E$, these lines are the proper transforms of the lines through
$\pi(E)$ and each of the images quadruple points, and the proper
transform of the conic through $\pi(E)$ and the image of the four
quadruple points.)  The linear system $\L_\P$ must thus satisfy the
following conditions:
\begin{enumerate}
  \item $\L_\P$ consists of degree $17$ curves
  \item $\L_\P$ contains two general $7$-tuple points, four general $6$-tuple points, and a simple point
  \item $\L_\P$ contains five triple points on a line.
\end{enumerate}
The last condition is a result of Lemma~\ref{le:match}.

We now must show that $\L_\P$, characterized above, must be empty.  To
do this, we degenerate the plane $\P$ into a reducible surface
$\P'\cup\F'$ just as before, this time letting the new component $\F'$
contain one of the triple points on the line, the two general $7$-tuple
points, and one $6$-tuple point (see Figure~\ref{steps}), and setting
$a=14$.  $\L_\P'$ is then a degree $14$ linear system with three
general $6$-tuple points, a general simple point, and eight points of
multiplicity 3 and 4 distributed on two lines according to the second
picture in Figure~1.

We repeat this twice more, according to the figure below.  At this
point we are left with a degree $8$ linear system, with ten points
distributed onto four lines and one general simple point according to
the forth diagram in the figure.  We specialize this simple point onto
the line through the two simple and two triple points; this line
splits off exactly once leaving us with $\L_7(2^4,3^4)$.  This must
contain in its base locus the $(-1)$-curve of degree $4$ through three
double points and five simple points.  Residual to his curve we have
$\L_3(1^7,2)$ which is clearly empty.

\begin{figure}[ht]
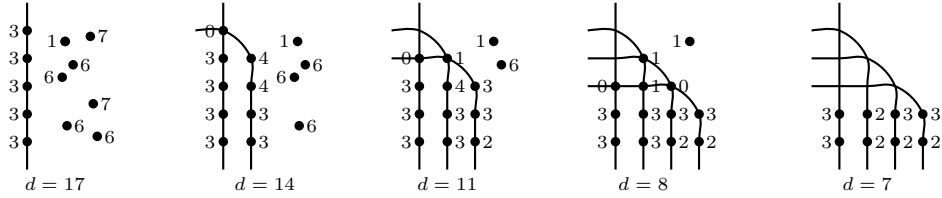
\tiny\label{steps}
\psset{labelsep=.10cm,unit=.37cm}
\pspicture(0,0)(4,6)
\pscurve(1,0)(1,3)(1,6)
\psdots(1,1)(1,2)(1,3)(1,4)(1,5)
\uput[180](1,1){3}
\uput[180](1,2){3}
\uput[180](1,3){3}
\uput[180](1,4){3}
\uput[180](1,5){3}
\psdots(3.38,2.37)(3.27,4.79)(3.52,1.19)(2.37,4.61)(2.26,3.32)(2.43,1.58)(2.65,3.77)(2.37,4.61)
\uput[0](3.38,2.37){7}
\uput[0](3.27,4.79){7}
\uput[0](3.52,1.19){6}
\uput[180](2.26,3.32){6}
\uput[0](2.43,1.58){6}
\uput[0](2.65,3.77){6}
\uput[180](2.37,4.61){1}
\uput[270](2,0){$d=17$}
\endpspicture
\hfill
\pspicture(0,0)(4,6)
\pscurve(1,0)(1,3)(1,6)
\pscurve(2,0)(2,1)(2,2)(2,3)(2,4)(1,5)(0,5)
\psdots(1,1)(1,2)(1,3)(1,4)(1,5)(2,1)(2,2)(2,3)(2,4)
\uput[180](1,1){3}
\uput[180](1,2){3}
\uput[180](1,3){3}
\uput[180](1,4){3}
\uput[180](1,5){0}
\uput[0](2,1){3}
\uput[0](2,2){3}
\uput[0](2,3){4}
\uput[0](2,4){4}
\psdots(3.56,3.32)(3.73,1.58)(3.95,3.77)(3.67,4.61)
\uput[180](3.56,3.32){6}
\uput[0](3.73,1.58){6}
\uput[0](3.95,3.77){6}
\uput[180](3.67,4.61){1}
\uput[270](2.5,0){$d=14$}
\endpspicture
\hfill
\pspicture(0,0)(4,6)
\pscurve(1,0)(1,3)(1,6)
\pscurve(2,0)(2,1)(2,2)(2,3)(2,4)(1,5)(0,5)
\pscurve(3,0)(3,1)(3,2)(3,3)(2,4)(1,4)(0,4)
\psdots(1,1)(1,2)(1,3)(1,4)(2,1)(2,2)(2,3)(2,4)(3,1)(3,2)(3,3)
\uput[180](1,1){3}
\uput[180](1,2){3}
\uput[180](1,3){3}
\uput[180](1,4){0}
\uput[0](2,1){3}
\uput[0](2,2){3}
\uput[0](2,3){4}
\uput[0](2,4){1}
\uput[0](3,1){2}
\uput[0](3,2){3}
\uput[0](3,3){3}
\psdots(3.95,3.77)(3.67,4.61)
\uput[0](3.95,3.77){6}
\uput[180](3.67,4.61){1}
\uput[270](2,0){$d=11$}
\endpspicture
\hfill
\pspicture(0,0)(5,6)
\psdots(1,1)(1,2)(1,3)(2,1)(2,2)(2,3)(2,4)(3,1)(3,2)(3,3)(4,1)(4,2)
\pscurve(1,0)(1,1)(1,2)(1,3)(1,4)(1,5)(1,6)
\pscurve(2,0)(2,1)(2,2)(2,3)(2,4)(1,5)(0,5)
\pscurve(3,0)(3,1)(3,2)(3,3)(2,4)(1,4)(0,4)
\pscurve(4,0)(4,1)(4,2)(3,3)(2,3)(1,3)(0,3)
\uput[180](1,1){3}
\uput[180](1,2){3}
\uput[180](1,3){0}
\uput[0](2,1){3}
\uput[0](2,2){3}
\uput[0](2,3){1}
\uput[0](2,4){1}
\uput[0](3,1){2}
\uput[0](3,2){3}
\uput[0](3,3){0}
\uput[0](4,1){2}
\uput[0](4,2){3}
\psdots(3.67,4.61)
\uput[180](3.67,4.61){1}
\uput[270](2,0){$d=8$}
\endpspicture
\hfill
\pspicture(0,0)(5,6)
\pscurve(1,0)(1,1)(1,2)(1,3)(1,4)(1,5)(1,6)
\pscurve(2,0)(2,1)(2,2)(2,3)(2,4)(1,5)(0,5)
\pscurve(3,0)(3,1)(3,2)(3,3)(2,4)(1,4)(0,4)
\pscurve(4,0)(4,1)(4,2)(3,3)(2,3)(1,3)(0,3)
\psdots(1,1)(1,2)(2,1)(2,2)(3,1)(3,2)(4,1)(4,2)
\uput[180](1,1){3}
\uput[180](1,2){3}
\uput[0](2,1){2}
\uput[0](2,2){2}
\uput[0](3,1){2}
\uput[0](3,2){3}
\uput[0](4,1){2}
\uput[0](4,2){3}
\uput[270](2,0){$d=7$}
\endpspicture
\caption{$\L_{21}(1,6^4,7^6)$ is empty}
\end{figure}

The case $\L_{22}(2, 6, 7^{9})$ is only slightly different and is
outlined in Figure~2.  This time we degenerate the plane five times,
using the nine $7$-tuple points and the $6$-tuple point.  The final
result ($d=7$) has fourteen points, mostly simple, distributed among
five lines.  In the last step, we specialize two triple points and one
double point onto a line, after which more are forced to appear in the
base locus than the degree of the residual linear system.  The
original linear system must be empty.  The details are not difficult
and are left to the reader.

\begin{figure}[ht]
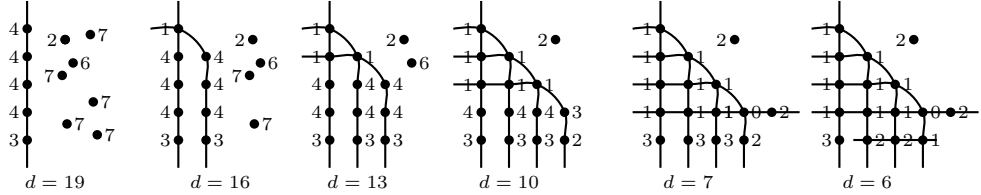
\tiny\label{steps22}
\psset{labelsep=.1cm,unit=.37cm}
\pspicture(0,0)(4,6)
\pscurve(1,0)(1,3)(1,6)
\psdots(1,1)(1,2)(1,3)(1,4)(1,5)
\uput[180](1,1){3}
\uput[180](1,2){4}
\uput[180](1,3){4}
\uput[180](1,4){4}
\uput[180](1,5){4}
\psdots(3.38,2.37)(3.27,4.79)(3.52,1.19)(2.37,4.61)(2.26,3.32)(2.43,1.58)(2.65,3.77)(2.37,4.61)
\uput[0](3.38,2.37){7}
\uput[0](3.27,4.79){7}
\uput[0](3.52,1.19){7}
\uput[180](2.26,3.32){7}
\uput[0](2.43,1.58){7}
\uput[0](2.65,3.77){6}
\uput[180](2.37,4.61){2}
\uput[270](2,0){$d=19$}
\endpspicture
\hfill
\pspicture(0,0)(4,6)
\pscurve(1,0)(1,3)(1,6)
\pscurve(2,0)(2,1)(2,2)(2,3)(2,4)(1,5)(0,5)
\psdots(1,1)(1,2)(1,3)(1,4)(1,5)(2,1)(2,2)(2,3)(2,4)
\uput[180](1,1){3}
\uput[180](1,2){4}
\uput[180](1,3){4}
\uput[180](1,4){4}
\uput[180](1,5){1}
\uput[0](2,1){3}
\uput[0](2,2){4}
\uput[0](2,3){4}
\uput[0](2,4){4}
\psdots(3.56,3.32)(3.73,1.58)(3.95,3.77)(3.67,4.61)
\uput[180](3.56,3.32){7}
\uput[0](3.73,1.58){7}
\uput[0](3.95,3.77){6}
\uput[180](3.67,4.61){2}
\uput[270](2.5,0){$d=16$}
\endpspicture
\hfill
\pspicture(0,0)(4,6)
\pscurve(1,0)(1,3)(1,6)
\pscurve(2,0)(2,1)(2,2)(2,3)(2,4)(1,5)(0,5)
\pscurve(3,0)(3,1)(3,2)(3,3)(2,4)(1,4)(0,4)
\psdots(1,1)(1,2)(1,3)(1,4)(1,5)(2,1)(2,2)(2,3)(2,4)(3,1)(3,2)(3,3)
\uput[180](1,1){3}
\uput[180](1,2){4}
\uput[180](1,3){4}
\uput[180](1,4){1}
\uput[180](1,5){1}
\uput[0](2,1){3}
\uput[0](2,2){4}
\uput[0](2,3){4}
\uput[0](2,4){1}
\uput[0](3,1){3}
\uput[0](3,2){4}
\uput[0](3,3){4}
\psdots(3.95,3.77)(3.67,4.61)
\uput[0](3.95,3.77){6}
\uput[180](3.67,4.61){2}
\uput[270](2,0){$d=13$}
\endpspicture
\hfill
\pspicture(0,0)(5,6)
\psdots(1,1)(1,2)(1,3)(1,4)(1,5)(2,1)(2,2)(2,3)(2,4)(3,1)(3,2)(3,3)(4,1)(4,2)
\pscurve(1,0)(1,1)(1,2)(1,3)(1,4)(1,5)(1,6)
\pscurve(2,0)(2,1)(2,2)(2,3)(2,4)(1,5)(0,5)
\pscurve(3,0)(3,1)(3,2)(3,3)(2,4)(1,4)(0,4)
\pscurve(4,0)(4,1)(4,2)(3,3)(2,3)(1,3)(0,3)
\uput[180](1,1){3}
\uput[180](1,2){4}
\uput[180](1,3){1}
\uput[180](1,4){1}
\uput[180](1,5){1}
\uput[0](2,1){3}
\uput[0](2,2){4}
\uput[0](2,3){1}
\uput[0](2,4){1}
\uput[0](3,1){3}
\uput[0](3,2){4}
\uput[0](3,3){1}
\uput[0](4,1){2}
\uput[0](4,2){3}
\psdots(3.67,4.61)
\uput[180](3.67,4.61){2}
\uput[270](2,0){$d=10$}
\endpspicture
\hfill
\pspicture(0,0)(5,6)
\pscurve(1,0)(1,1)(1,2)(1,3)(1,4)(1,5)(1,6)
\pscurve(2,0)(2,1)(2,2)(2,3)(2,4)(1,5)(0,5)
\pscurve(3,0)(3,1)(3,2)(3,3)(2,4)(1,4)(0,4)
\pscurve(4,0)(4,1)(4,2)(3,3)(2,3)(1,3)(0,3)
\pscurve(0,2)(3,2)(6,2)
\psdots(1,1)(1,2)(1,3)(1,4)(1,5)(2,1)(2,2)(2,3)(2,4)(3,1)(3,2)(3,3)(4,1)(4,2)(5,2)
\uput[180](1,1){3}
\uput[180](1,2){1}
\uput[180](1,3){1}
\uput[180](1,4){1}
\uput[180](1,5){1}
\uput[0](2,1){3}
\uput[0](2,2){1}
\uput[0](2,3){1}
\uput[0](2,4){1}
\uput[0](3,1){3}
\uput[0](3,2){1}
\uput[0](3,2){1}
\uput[0](3,3){1}
\uput[0](4,1){2}
\uput[0](4,2){0}
\uput[0](5,2){2}
\psdots(3.67,4.61)
\uput[180](3.67,4.61){2}
\uput[270](2,0){$d=7$}
\endpspicture
\hfill
\pspicture(0,0)(5,6)
\pscurve(1,0)(1,1)(1,2)(1,3)(1,4)(1,5)(1,6)
\pscurve(2,0)(2,1)(2,2)(2,3)(2,4)(1,5)(0,5)
\pscurve(3,0)(3,1)(3,2)(3,3)(2,4)(1,4)(0,4)
\pscurve(4,0)(4,1)(4,2)(3,3)(2,3)(1,3)(0,3)
\pscurve(0,2)(3,2)(6,2)
\pscurve(1.5,1)(3,1)(4.5,1)
\psdots(1,1)(1,2)(1,3)(1,4)(1,5)(2,1)(2,2)(2,3)(2,4)(3,1)(3,2)(3,3)(4,1)(4,2)(5,2)
\uput[180](1,1){3}
\uput[180](1,2){1}
\uput[180](1,3){1}
\uput[180](1,4){1}
\uput[180](1,5){1}
\uput[0](2,1){2}
\uput[0](2,2){1}
\uput[0](2,3){1}
\uput[0](2,4){1}
\uput[0](3,1){2}
\uput[0](3,2){1}
\uput[0](3,3){1}
\uput[0](4,1){1}
\uput[0](4,2){0}
\uput[0](5,2){2}
\psdots(3.67,4.61)
\uput[180](3.67,4.61){2}
\uput[270](2,0){$d=6$}
\endpspicture
\caption{$\L_{22}(2,6,7^9)$ is empty}
\end{figure}

\end{proof}

\begin{lemma}\label{le:save}
$\L_{22}(2,6^{13})$ is empty.
\end{lemma}
\begin{proof}
Let $\pi:X\to\P^2$ be the blow-up of the projective plane at the
thirteen $6$-tuple points.  By the triangular checker game, we know that
$\L_{22}(6^{13})$ is non-special and so the linear system
$|22H-\sum_{i=1}^{13} 6E_i|$ on $X$ gives rise to a map
$\phi:X\to\P^2$.

Assume for a contradiction that $\L=\L_{22}(2,6^{13})$ is non-empty.
The differential of $\phi$ fails to be injective at a general point in
$\P^2$, thus the image $f(X)\in\P^2$ is a curve, say of degree $d$.
Since $\L=f^*\O_{\P^2}(1)$, the general element of $\L_{22}(2,6^{13})$
has $d$ components of the same degree.  The degree $d$ must be a
divisor of $22$, and a case-by-case analysis shows that this is
impossible.  Clearly $d$ cannot be $1$ or $22$.  If $d=11$, then
$\L_{22}(6^{13})$ is comprised of $11$ conics, which is impossible.
If $d=2$, then $\L_{22}(6^{13})$ is comprised of two curves in
$\L_{11}(1,3^{13})$, which is empty by previous calculations.
\end{proof}

These last two lemmas suffice to prove that the remaining linear
systems in Table~\ref{42} are non-special and empty.

\section*{References}

\bibliographystyle{amsxport}

\begin{biblist} 
\bib{MR2129797}{article}{
    author={Bocci, C.},
    author={Miranda, R.},
     title={Topics on interpolation problems in algebraic geometry},
   journal={Rend. Sem. Mat. Univ. Politec. Torino},
    volume={62},
      date={2004},
    number={4},
     pages={279\ndash 334},
      issn={0373-1243},
    review={MR2129797},
}
\bib{MR2000m:14005}{article}{
  author={Ciliberto, Ciro},
  author={Miranda, Rick},
  title={Degenerations of planar linear systems},
  journal={J. Reine Angew. Math.},
   volume={501},
  date={1998},
  pages={191\ndash 220},
  issn={0075-4102},
  review={MR2000m:14005},
}
\bib{MR2000m:14006}{article}{
  author={Ciliberto, Ciro},
  author={Miranda, Rick},
  title={Linear systems of plane curves with base points of equal multiplicity},
  journal={Trans. Amer. Math. Soc.},
  volume={352},
  date={2000},
  number={9},
  pages={4037\ndash 4050},
  issn={0002-9947},
  review={MR 2000m:14006},
}
\bib{matching}{article}{
  author={Ciliberto, Ciro},
  author={Miranda, Rick},
  title={Matching Conditions for Degenerating Plane Curves and Applications},
  note={in preparation},
}
\bib{CCMO}{article}{
  author={Ciliberto, Ciro},
  author={Cioffi, Francesca},
  author={Miranda, Rick},
  author={Orecchia, Ferruccio},
  title={Bivariate Hermite interpolation via computer algebra and algebraic geometry techniques},
  journal={Pubblicazioni del Dipartimento di Matematica e Applicazioni di Napoli},
  number={26},
  year={2002},
}
\bib{MR91a:14007}{article}{
  author={Gimigliano, Alessandro},
  title={Our thin knowledge of fat points},
  booktitle={The Curves Seminar at Queen's, Vol.\ VI (Kingston, ON, 1989)},
  series={Queen's Papers in Pure and Appl. Math.},
  volume={83},
  pages={Exp.\ No.\ B, 50},
  publisher={Queen's Univ.},
  place={Kingston, ON},
  date={1989},
  review={MR 91a:14007},
}
\bib{MR846019}{article}{
  author={Harbourne, Brian},
  title={The geometry of rational surfaces and Hilbert functions of points in the plane},
  booktitle={Proceedings of the 1984 Vancouver conference in algebraic geometry},
  series={CMS Conf. Proc.},
  volume={6},
  pages={95\ndash 111},
  publisher={Amer. Math. Soc.},
  place={Providence, RI},
  date={1986},
  review={MR846019 (87k:14041)},
}
\bib{MR2003f:13032}{article}{
  author={Harbourne, Brian},
  title={Problems and progress: a survey on fat points in $\mathbb {P}^2$},
  booktitle={Zero-dimensional schemes and applications (Naples, 2000)},
  series={Queen's Papers in Pure and Appl. Math.},
  volume={123},
  pages={85\ndash 132},
  publisher={Queen's Univ.},
  place={Kingston, ON},
  date={2002},
  review={MR 2003f:13032},
}
\bib{MR932136}{article}{
    author={Hirschowitz, Andr{\'e}},
     title={Existence de faisceaux r\'eflexifs de rang deux sur ${\bf P}\sp
            3$ \`a bonne cohomologie},
  language={French},
   journal={Inst. Hautes \'Etudes Sci. Publ. Math.},
    number={66},
      date={1988},
     pages={105\ndash 137},
      issn={0073-8301},
    review={MR932136 (89c:14019)},
}
\bib{MR993223}{article}{
  author={Hirschowitz, Andr{\'e}},
  title={Une conjecture pour la cohomologie des diviseurs sur les surfaces rationnelles g\'en\'eriques},
  language={French},
  journal={J. Reine Angew. Math.},
  volume={397},
  date={1989},
  pages={208\ndash 213},
  issn={0075-4102},
  review={MR993223 (90g:14021)},
}
\bib{MR21:4151}{article}{
  author={Nagata, Masayoshi},
  title={On the $14$-th problem of Hilbert},
  journal={Amer.~J.~Math.},
  volume={81},
  date={1959},
  pages={766\ndash 772},
  review={MR 21 \#4151},
}
\end{biblist}

\end{document}